\documentclass[11pt]{article}
\usepackage{mathrsfs}
\usepackage{}
\usepackage{amsfonts}
\usepackage{amssymb}
\usepackage{amsmath}

\def\sqr#1#2{{\vcenter{\vbox{\hrule height.#2pt
              \hbox{\vrule width.#2pt height#1pt \kern#1pt \vrule
width.#2pt}
              \hrule height.#2pt}}}}
\def\signed #1{{\unskip\nobreak\hfil\penalty50
              \hskip2em\hbox{}\nobreak\hfil#1
              \parfillskip=0pt \finalhyphendemerits=0 \par}}
\def\endpf{\signed {$\sqr69$}}
\def\3n{\negthinspace \negthinspace \negthinspace }
\def\2n{\negthinspace \negthinspace }
\def\1n{\negthinspace }

\def\O{\Omega}
\def\({\Big (}
\def\){\Big )}
\def\[{\Big[}
\def\]{\Big]}

\def\be{\begin{equation}}
\def\bel{\begin{equation}\label}
\def\ee{\end{equation}}
\def\bea{\begin{eqnarray}}
\def\eea{\end{eqnarray}}
\def\bt{\begin{theorem}}
\def\et{\end{theorem}}
\def\bc{\begin{corollary}}
\def\ec{\end{corollary}}
\def\bl{\begin{lemma}}
\def\el{\end{lemma}}
\def\bp{\begin{proposition}}
\def\ep{\end{proposition}}
\def\br{\begin{remark}}
\def\er{\end{remark}}
\def\ba{\begin{array}}
\def\ea{\end{array}}
\def\bd{\begin{definition}}
\def\ed{\end{definition}}

\newtheorem{lemma}{Lemma}[section]
\newtheorem{remark}{Remark}[section]

\newtheorem{theorem}{Theorem}[section]
\newtheorem{corollary}{Corollary}[section]

\newtheorem{definition}{Definition}[section]
\newtheorem{proposition}{Proposition}[section]

\makeatletter
   
   \@addtoreset{equation}{section}
\makeatother
\allowdisplaybreaks
\newcommand{\upcite}[1]{\textsuperscript{\cite{#1}}}
\begin{document}

\title{\bf Hierarchic control for the coupled fourth order parabolic equations }

\author{Yating Li  and   Muming Zhang*}

\date{School of Mathematics and Statistics, Northeast Normal
University, \\
Changchun 130024, China.}

\maketitle

\renewcommand{\thefootnote}{\fnsymbol{footnote}}
\footnote[0]{ *Corresponding author (M. Zhang): zhangmm352@nenu.edu.cn.}
\footnote[0]{This work is supported by National Key R\&D Program of China under grant 2023YFA1009002, and NSF of China under grants 12001094, 12001087 and 11971179, and Fundamental Research Funds for the Central Universities under grant 2412020QD027.  }

\begin{abstract}

In this paper, we obtain a null controllability result for a
 coupled fourth order parabolic system based on the Stackelberg-Nash strategies. For this purpose, we first prove the existence and uniqueness of Nash equilibrium pair of the original system  and its explicit expression is provided. Next, we investigate the null controllability of Nash equilibrium to the corresponding optimal system.  By duality theory, we establish an observability estimate for the coupled fourth order parabolic system. Such an estimate is obtained by a new global Carleman estimate we derived.
\end{abstract}


\section{Introduction and main result}

Let $T>0$ and $N\geq2$. Assume that $\Omega\subset\mathbb{R}^{N}$ is a bounded domain whose boundary
$\partial\Omega$ is regular enough. Let $\omega,\omega_{1}$  and $\omega_{2}$ be three given nonempty open subsets of $\Omega$, and $\omega_{i}\cap\omega=\emptyset  (i=1,2).$ Denote by $\chi_{\omega}$ the characteristic function of $\omega$.  Put $Q=\Omega\times(0,T)$ and $\Sigma=\partial\Omega\times(0,T)$. Unless otherwise stated, we assume that the set of indicators $i$ mentioned in this paper contains $i=1$ and $i=2$.

Consider the following coupled fourth order parabolic equations:
\begin{equation}\label{1.1}
\left\{\begin{array}{ll}
y_{1,t}+\Delta^{2}y_{1}+a_{11}y_{1}+a_{12}y_{2}=g\chi_{\omega}+h_{1}\chi_{\omega_{1}}+h_{2}\chi_{\omega_{2}}& \mbox{ in }~Q,\\[3mm]
y_{2,t}+\Delta^{2}y_{2}+a_{21}y_{1}+a_{22}y_{2}=0& \mbox{ in }~Q,\\[3mm]
y_{j}=\Delta y_{j}=0&\mbox{ on }~\Sigma,~j=1,2,\\[3mm]
y_{j}(x,0)=y_{j}^{0}(x)&\mbox{ in }~\Omega,~j=1,2,
\end{array}\right.
\end{equation}
\\where $g\in L^{2}(\omega\times(0,T))$ is the leader control, $h_{i}\in L^{2}(\omega_{i}\times(0,T))~(i=1,2)$ are the follower controls,
$y=y(\cdot,\cdot;g,h_{1},h_{2})=(y_{1},y_{2})^{T}$ is the state.
 $a_{ij}=a_{ij}(\cdot,\cdot)\in L^{\infty}(Q)~(i,j=1,2)$, $y_{j}^{0}\in L^{2}(\Omega)(j=1,2)$ are given functions, and we write $y^{0}=(y_{1}^{0},y_{2}^{0})$.

\medskip
The aim of this paper is to study the hierarchical control of the coupled fourth order parabolic equations (\ref{1.1}). For this purpose,  we need to introduce two types of  functionals.
\\(1)~The main cost functional:
\begin{equation}\label{1.2}
 J(g)=\frac{1}{2}\int_{\omega\times(0,T)}|g|^{2}dxdt.
\end{equation}
(2)~The secondary cost functionals:
\begin{equation}\label{1.3}
 J_{i}(g~;h_{1},h_{2})=\frac{\alpha_{i}}{2}\!\int_{O_{i,d}\times(0,T)}\!|y_{1}-y_{1,d}^{i}|^{2}+|y_{2}-y_{2,d}^{i}|^{2}dxdt+\frac{\mu_{i}}{2}\!\int_{\omega_{i}\times(0,T)}\!\rho_{*}^{2}|h_{i}|^{2}dxdt,
\end{equation}
where $\alpha_{i}$ and $\mu_{i}$ are positive constants, $\rho_{*}=\rho_{*}(t)$  is given in (\ref{z3.6}) as a positive function, $O_{i,d}\subseteq\Omega (i=1,2)$ are nonempty open subsets of $\O$, which represent  two observation domains, $y_{d}^{i}=(y_{1,d}^{i},y_{2,d}^{i})\in L^{2}(O_{i,d}\times(0,T))^{2}$ are given functions, and $y=y(\cdot,\cdot;g,h_{1},h_{2})$ is the solution of system (\ref{1.1}) corresponding to the leader $g$ and the followers $(h_{1},h_{2})$.
\medskip

We give the  definition of Nash equilibrium pair for the secondary cost functionals $J_{i}$ as follows.

\begin{definition}\label{d3.1}
For any given leader control g,  a follower control pair $(\bar{h}_{1},\bar{h}_{2})$ is said to be a Nash equilibrium for $(J_{1},J_{2})$, if the following identities
\begin{equation}\label{1.4}
J_{1}(g;\bar{h}_{1},\bar{h}_{2})=\min_{h_{1}} J_{1}(g;{h}_{1},\bar{h}_{2}),~~\forall h_{1}\in L^{2}(\omega_{1}\times(0,T)),
\end{equation}
\begin{equation}\label{1.5}
J_{2}(g;\bar{h}_{1},\bar{h}_{2})=\min_{h_{2}} J_{2}(g;\bar{h}_{1},{h}_{2}),~~\forall h_{2}\in L^{2}(\omega_{2}\times(0,T))
\end{equation}
hold.
\end{definition}

The main interests in this paper center on the null controllability of coupled fourth order parabolic system (\ref{1.1}) with Stackelberg-Nash strategies, that is, we have two goals:
\\1. For any given leader control $g$, prove that the Nash equilibrium for $(J_1,J_2)$ exists, denoted by $(\bar{h}_{1}(g),\bar{h}_{2}(g))$.  

\noindent
2.  Prove the existence of  a leader control $\bar{g}\in L^{2}(\omega\times(0,T))$  that satisfies
\begin{equation}\label{1.6}
J(\bar{g})=\min_{g} J(g;\bar{h}_{1}(g),\bar{h}_{2}(g)),~~\forall g\in L^{2}(\omega\times(0,T)),
\end{equation}
subjecting to the restriction 
\begin{equation}\label{1.7}
y(\cdot,T;\bar{g},\bar{h}_{1}(\bar{g}),\bar{h}_{2}(\bar{g}))=0~\mbox{ in }~\Omega,
\end{equation}
where $y=y(\cdot,\cdot;\bar{g},\bar{h}_{1}(\bar{g}),\bar{h}_{2}(\bar{g}))$ is the solution of system (\ref{1.1}) corresponding to the leader $\bar{g}$ and its  Nash equilibrium  $(\bar{h}_{1}(\bar{g}),\bar{h}_{2}(\bar{g}))$.

\medskip
Our main result is stated as follows.

\begin{theorem}\label{t1}
Assume that $O_{1,d}=O_{2,d}$, denoted as $O_{d}$, $O_{d}\cap\omega\neq\emptyset$ and $\mu_{i}$ (i=1,2) are sufficiently large. If
$$a_{21}\geq a_{0}>0\ \mbox{or}\ -a_{21}\geq a_{0}>0\ \mbox{in}\ (O_{d}\cap\omega)\times(0,T)$$
and
\begin{equation}\label{1.8}
\int_{O_{d}\times(0,T)}|y_{j,d}^{i}|^{2}dxdt<+\infty,\ i,j=1,2,
\end{equation}
then for any $y^{0}=(y_{1}^{0},y_{2}^{0})\in L^{2}(\Omega)^{2}$, there exists a control $\bar{g}\in L^{2}(\omega\times(0,T))$ and the corresponding Nash equilibrium $(\bar{h}_{1}(\bar{g}),\bar{h}_{2}(\bar{g}))$ such that the solution of (\ref{1.1}) satisfies (\ref{1.6}) and (\ref{1.7}).
\end{theorem}
\medskip

In classical control theory, there are two fundamental problems: controllability and optimal control.  Roughly speaking, controllability means that we can find a control to achieve a predetermined goal, while optimal control means finding the best way to achieve a goal in a sense. However, controllability and optimal control can be involved simultaneously in practical problems, i.e., multi-objective control problems.
 It is necessary to exert more than one control in order to accomplish multiple objectives. Therefore, hierarchical control is proposed by Lions in \cite{CO} to solve the control problem of multiple objectives. These objectives may be cooperative or noncooperative, and achieving multiple objectives means finding some sense of equilibrium. Here we refer to the Nash equilibrium in Game theory, which was proposed by J. F. Nash in \cite{JN}, also called the noncooperative game equilibrium.
The main idea of the Stackelberg-Nash strategy (see \cite{stack}) is to choose its own strategy according to the  possible strategies of other participants to ensure that one can obtain the maximum benefit. 
The leader makes decisions first, and followers choose strategies based on leader's decisions. 
There are numerous studies on the Stackelberg-Nash controllability for parabolic equations (see \cite{afg2017},\cite{afs2015},\cite{afs2020} and the rich references cited therein).  In \cite{jde2019} and \cite{jde20231-63}, the Stackelberg-Nash  exact controllability for the Kuramoto-Sivashinsky equations,  a special class of nonlinear  fourth-order parabolic equations, was discussed. It is worth mentioning that, \cite{2022arxiv} studied the hierarchical exact controllability problem of the fourth-order linear and semilinear parabolic equations. About coupled parabolic systems, in \cite{htp2016}, the authors analyzed the Stackelberg-Nash null controllability of a coupled linear parabolic system. Further, the similar problem for linear and semilinear degenerate parabolic equations was investigated in \cite{aaf2018}.  However, to our best knowledge, there are no references considering the Stackelberg-Nash null controllability of system (\ref{1.1}).

In this paper, we study the hierarchic control based on Stackelberg-Nash strategy for the coupled fourth order parabolic equations.  Hierarchical control problems are generally solved by first proving the existence of Nash equilibrium pair of the original system, then proving the null controllability of Nash equilibrium to the corresponding optimal system. Up to now, several authors have already studied the problem of controllability for the fourth order parabolic equations of one and higher dimensions (see, for instance, \cite{2016carreno},\cite{2006chouguo},\cite{2002guo},\cite{2020kassab},\cite{2009yu},\cite{2012zhou}). We refer to \cite{2019guerrero} for the first result for the Carleman estimate of a fourth order parabolic equation in dimension $n\geq 2$. Further, the controllability of the coupled fourth order parabolic equations was investigated in \cite{2019carrenocma}. Recently, there are some authors discussing the controllability of fourth order stochastic parabolic equations (see \cite{lwang}).

The rest of this paper is organized as follows. In Section 2, we shall prove the existence and uniqueness of a Nash equilibrium for $(J_1,J_2)$.
In Section 3, we transform the existence of leader control satisfying (\ref{1.6}) and (\ref{1.7}) into an observability estimate for a coupled parabolic system.  Section 4 is devoted to deriving a new Carleman estimate for coupled  fourth order parabolic equations. Based on this, the
observability estimate is established in Section 5.

\section{Nash equilibrium}
In this section, we establish the existence and uniqueness of the Nash equilibrium and give its explicit expression.

\subsection{Existence and uniqueness of a Nash equilibrium}\label{z1}

Define the spaces
$$H_{i}=L^{2}(\omega_{i}\times(0,T)),~i=1,2,~~~~H=H_{1}\times H_{2}.$$
and put the operators~$\Lambda_{i}:H_{i}\rightarrow L^{2}(Q)^{2}$ by $\Lambda_{i}h_{i}=y^{i}$, where~$y^{i}=(y_{1}^{i},y_{2}^{i})~(i=1,2)$~are the solutions of the following systems
\begin{equation}\label{z3.1}
\left\{\begin{array}{ll}
y_{1,t}^{i}+\Delta^{2}y_{1}^{i}+a_{11}y_{1}^{i}+a_{12}y_{2}^{i}=h_{i}\chi_{\omega_{i}}& \mbox{ in }~Q,\\[3mm]
y_{2,t}^{i}+\Delta^{2}y_{2}^{i}+a_{21}y_{1}^{i}+a_{22}y_{2}^{i}=0& \mbox{ in }~Q,\\[3mm]
y_{j}^{i}=\Delta y_{j}^{i}=0&\mbox{ on }~\Sigma,~j=1,2,\\[3mm]
y_{j}^{i}(x,0)=0&\mbox{ in }~\Omega,~j=1,2.
\end{array}\right.
\end{equation}
For any $g\in L^{2}(\omega\times(0,T))$, we write the solution of (\ref{1.1}) as follows
$$y=y^{1}+y^{2}+u(g)=\Lambda_{1}h_{1}+\Lambda_{2}h_{2}+u(g),$$
where $u(g)=(u_{1}(g),u_{2}(g))^{T}$ is the solution of
\begin{equation}\label{z3.2}
\left\{\begin{array}{ll}
u_{1,t}+\Delta^{2}u_{1}+a_{11}u_{1}+a_{12}u_{2}=g\chi_{\omega}& \mbox{ in }~Q,\\[3mm]
u_{2,t}+\Delta^{2}u_{2}+a_{21}u_{1}+a_{22}u_{2}=0& \mbox{ in }~Q,\\[3mm]
u_{j}=\Delta u_{j}=0&\mbox{ on }~\Sigma,~j=1,2,\\[3mm]
u_{j}(x,0)=y_{j}^{0}(x)&\mbox{ in }~\Omega,~j=1,2.
\end{array}\right.
\end{equation}
Thus, the functionals (\ref{1.3}) can be written as
\begin{eqnarray*}
\begin{array}{rl}
&J_{i}(g~;h_{1},h_{2})=\displaystyle\frac{\alpha_{i}}{2}\!\!\int_{O_{i,d}\times(0,T)}\!\!\|\Lambda_{1}h_{1}\!+\!\Lambda_{2}h_{2}\!+\!u(g)\!-\!y_{d}^{i}\|^{2}dxdt+\displaystyle\frac{\mu_{i}}{2}\!\!\int_{\omega_{i}\times(0,T)}\!\rho_{*}^{2}|h_{i}|^{2}dxdt
\\[3mm]&\quad\quad\quad\quad\quad~=\displaystyle\frac{\alpha_{i}}{2}\!\!\int_{O_{i,d}\times(0,T)}\|\Lambda_{1}h_{1}\!+\!\Lambda_{2}h_{2}\!-\!\tilde y_{d}^{i}\|^{2}dxdt+\displaystyle\frac{\mu_{i}}{2}\!\!\int_{\omega_{i}\times(0,T)}\rho_{*}^{2}|h_{i}|^{2}dxdt,
\end{array}
\end{eqnarray*}
where $\displaystyle\tilde y_{d}^{i}=y_{d}^{i}-u(g)\chi_{O_{i,d}}~(i=1,2)$, and $\|\cdot\|$ stands for the usual Euclidian norm in~$\mathbb{R}^{2}$.

Since the functionals~$J_{1},J_{2}$~are strictly convex, by Definition \ref{d3.1}, it follows that ~$(\bar{h}_{1},\bar{h}_{2})$~is a Nash equilibrium for~$(J_{1},J_{2})$~if and only if
\begin{equation*}
\left(~\frac{\partial J_{1}}{\partial h_{1}}(g,\bar{h}_{1},\bar{h}_{2}),~h_{1}\right)=0,~\forall h_{1}\in H_{1};
~~\left(~\frac{\partial J_{2}}{\partial h_{2}}(g,\bar{h}_{1},\bar{h}_{2}),~h_{2}\right)=0,~\forall h_{2}\in H_{2},
\end{equation*}
where~$\left(\displaystyle~\frac{\partial J_{1}}{\partial h_{1}}(g,\bar{h}_{1},\bar{h}_{2}),~h_{1}\right)$~represents~$G\hat{a}teaux$~derivative of~$J_{1}$~at~$(g,\bar{h}_{1},\bar{h}_{2})$~along~$h_{1}$.
 That is
$$\displaystyle\mu_{1}\int_{\omega_{1}\times(0,T)}\rho_{*}^{2}\bar{h}_{1}h_{1}~dxdt+\alpha_{1}\int_{O_{1,d}\times(0,T)}(\Lambda_{1}\bar{h}_{1}+\Lambda_{2}\bar{h}_{2}-\tilde{y}_{d}^{1})\cdot\Lambda_{1}h_{1}~dxdt=0,~\forall h_{1}\in H_{1}.$$
Therefore,~$(\bar{h}_{1},\bar{h}_{2})$~is a Nash equilibrium if and only if
\begin{equation}\label{z3.5}
\mu_{i}\int_{\omega_{i}\times(0,T)}\rho_{*}^{2}\bar{h}_{i}h_{i}~dxdt+\alpha_{i}\int_{O_{i,d}\times(0,T)}(\Lambda_{1}\bar{h}_{1}+\Lambda_{2}\bar{h}_{2}-\tilde{y}_{d}^{i})\cdot\Lambda_{i}h_{i}~dxdt=0, ~i=1,2.
\end{equation}
For any~$(h_{1},h_{2})\in H $, it holds that 
\begin{equation*}
\displaystyle\mu_{i}(\rho_{*}^{2}\bar{h}_{i},h_{i})_{L^{2}(\omega_{i}\times(0,T))}+\alpha_{i}\Big{(}\Lambda_{i}^{*}\Big{[}(\Lambda_{1}\bar{h}_{1}+\Lambda_{2}\bar{h}_{2}){|}_{O_{i,d}}-\tilde{y}_{d}^{i}\Big{]},h_{i}\Big{)}_{L^{2}(\omega_{i}\times(0,T))}=0, ~i=1,2,
\end{equation*}
where~$(\cdot,\cdot)_{L^{2}(\Omega)}$~denotes the internal product in~$L^{2}(\Omega)$, and~$\Lambda_{i}^{*}\in \mathcal{L}(L^{2}(Q)^{2},H_{i})$~is the adjoint operator of~$\Lambda_{i}$. This implies
$$\mu_{i}\rho_{*}^{2}\bar{h}_{i}+\alpha_{i}\Lambda_{i}^{*}\Big{[}(\Lambda_{1}\bar{h}_{1}+\Lambda_{2}\bar{h}_{2}){|}_{O_{i,d}}\Big{]}=\alpha_{i}\Lambda_{i}^{*}\tilde{y}_{d}^{i},~i=1,2.$$ 
For any~$h=(h_{1},h_{2})$, we define the operator~$K=(K_{1},K_{2})$,  where $K_{1}\in \mathcal{L}(H,H_{1})$ and $K_{2}\in \mathcal{L}(H,H_{2})$ are given by 
$$K_{i}h=\mu_{i}\rho_{*}^{2}{h}_{i}+\alpha_{i}\Lambda_{i}^{*}\Big{[}(\Lambda_{1}{h}_{1}+\Lambda_{2}{h}_{2}){|}_{O_{i,d}}\Big{]}, ~i=1,2.$$
Therefore,~$\bar{h}=(\bar{h}_{1},\bar{h}_{2})$~is a Nash equilibrium pair equivalent to
\begin{equation}\label{z3.7}
K_{i}\bar h=\alpha_{i}\Lambda_{i}^{*}\tilde{y}_{d}^{i}, ~i=1,2.
\end{equation}
Then,
\begin{eqnarray}\label{z3.8}
\begin{array}{rl}
& (Kh,h)_{H}=(K_{1}h,h_{1})_{H_{1}}+(K_{2}h,h_{2})_{H_{2}}
\\[3mm]&\quad\quad\quad\quad=\displaystyle\sum_{i=1}^{2}\mu_{i}\|\rho_{*}h_{i}\|_{L^{2}(\omega_{i}\times(0,T))}^{2}+\alpha_{1}(\Lambda_{1}h_{1}+\Lambda_{2}h_{2},\Lambda_{1}h_{1})_{L^{2}(O_{1,d}\times(0,T))^{2}}
\\[3mm]&\quad\quad\quad\quad\quad+\alpha_{2}(\Lambda_{1}h_{1}+\Lambda_{2}h_{2},\Lambda_{2}h_{2})_{L^{2}(O_{2,d}\times(0,T))^{2}}.
\end{array}
\end{eqnarray}
By applying Cauchy inequality with~$\varepsilon$, we obtain
\begin{eqnarray}\label{z3.9}
\begin{array}{rl}
&\alpha_{1}(\Lambda_{1}h_{1}+\Lambda_{2}h_{2},\Lambda_{1}h_{1})_{L^{2}(O_{1,d}\times(0,T))^{2}}
\\[3mm]&=\alpha_{1}\|\Lambda_{1}h_{1}\|^{2}_{L^{2}(O_{1,d}\times(0,T))^{2}}+\alpha_{1}(\Lambda_{2}h_{2},\Lambda_{1}h_{1})_{L^{2}(O_{1,d}\times(0,T))^{2}}
\\[3mm]&\geq\alpha_{1}\|\Lambda_{1}h_{1}\|^{2}_{L^{2}(O_{1,d}\times(0,T))^{2}}-\displaystyle\frac{\alpha_{1}}{2\varepsilon}\|\Lambda_{1}h_{1}\|^{2}_{L^{2}(O_{1,d}\times(0,T))^{2}}-\displaystyle\frac{\alpha_{1}\varepsilon}{2}\|\Lambda_{2}h_{2}\|^{2}_{L^{2}(O_{1,d}\times(0,T))^{2}}
\\[3mm]&\geq-\displaystyle\frac{\alpha_{1}}{4}\|\Lambda_{2}h_{2}\|^{2}_{L^{2}(O_{1,d}\times(0,T))^{2}}\geq-\displaystyle\frac{\alpha_{1}}{4}\|\Lambda_{2}\chi_{O_{1,d}}\|^{2}_{H_{1,d}}\|h_{2}\|^{2}_{H_{2}},
\end{array}
\end{eqnarray}
where~$\varepsilon=1/2$, and~$\|\cdot\|_{H_{i,d}}$~denotes the norm in~$\displaystyle \mathcal{L}(H_{3-i},~L^{2}(O_{i,d}\times(0,T))^{2})$. We put~$\rho_{0}=\displaystyle\min_{t\in[0,T]}\rho_{*}(t)$. By~(\ref{z3.8})~and~(\ref{z3.9}),~we have
\begin{eqnarray*}
\begin{array}{rl}
&(Kh,h)_{H}\geq\mu_{1}\rho_{0}^{2}\|h_{1}\|_{H_{1}}^{2}+\mu_{2}\rho_{0}^{2}\|h_{2}\|_{H_{2}}^{2}-\displaystyle\frac{\alpha_{1}}{4}\|\Lambda_{2}\chi_{O_{1,d}}\|^{2}_{H_{1,d}}\|h_{2}\|^{2}_{H_{2}}
\\[3mm]&\quad\quad\quad\quad\quad-\displaystyle\frac{\alpha_{2}}{4}\|\Lambda_{1}\chi_{O_{2,d}}\|^{2}_{H_{2,d}}\|h_{1}\|^{2}_{H_{1}}.
\end{array}
\end{eqnarray*}
Then, for~$\mu_{1},\mu_{2}$ large enough satisfying
\begin{equation}\label{z3.11}
\mu_{1}>\frac{\alpha_{2}\|\Lambda_{1}\chi_{O_{2,d}}\|^{2}_{H_{2,d}}}{4\rho_{0}^{2}},~\mu_{2}>\frac{\alpha_{1}\|\Lambda_{2}\chi_{O_{1,d}}\|^{2}_{H_{1,d}}}{4\rho_{0}^{2}},
\end{equation}
we have
\begin{equation}\label{z3.12}
(Kh,h)_{H}\geq \tau\|h\|^{2}_{H},
\end{equation}
where~$\tau=\displaystyle\min_{i=1,2}\left\{\mu_{i}\rho_{0}^{2}-\frac{\alpha_{3-i}}{4}\|\Lambda_{i}\chi_{O_{3-i,d}}\|^{2}_{H_{3-i,d}}\right\}.$ By~(\ref{z3.11}),~$\tau>0$.

We introduce the functional~$F(h,q):H\times H\rightarrow \mathbb{R}$~by $F(h,q)=(Kh,q)_{H}.$~Then, according to the definition of~$K$,~$F$~is bounded and bilinear. Moreover, by (\ref{z3.12}), $F$~is coercive. Applying~Lax-Milgram~theorem, for any given~$v\in H$, there exists a unique ~$\bar{h}=(\bar{h}_{1},\bar{h}_{2})\in H$~such that
$$F(\bar{h},q)=(v,q)_{H},~~\forall q\in H.$$
Set $v=\displaystyle(\alpha_{1}\Lambda_{1}^{*}\tilde{y}_{d}^{1},\alpha_{2}\Lambda_{2}^{*}\tilde{y}_{d}^{2})^{T}$, then we have~$K\bar{h}=v$~and (\ref{z3.7})~holds, where the function~$\bar{h}$~ is the desired Nash equilibrium pair. Thus we deduce the existence and uniqueness of the Nash equilibrium related to~$(J_{1},J_{2})$.

\subsection{Explicit expression of the Nash equilibrium pair}

We have proved that, for~$\mu_{1} $, $\mu_{2}$~large enough, there exists a unique Nash equilibrium~$(\bar{h}_{1},\bar{h}_{2})$ for $(J_{1},J_{2})$ in Subsection \ref{z1}. Moreover, by~(\ref{z3.5}), $(\bar{h}_{1},\bar{h}_{2})$~is the Nash equilibrium pair equivalent to 
\begin{equation}\label{z3.13}
\mu_{i}\int_{\omega_{i}\times(0,T)}\rho_{*}^{2}\bar{h}_{i}\hat{h}_{i}dxdt+\alpha_{i}\int_{O_{i,d}\times(0,T)}(\Lambda_{1}\bar{h}_{1}+\Lambda_{2}\bar{h}_{2}-\tilde{y}_{d}^{i})\cdot\Lambda_{i}\hat{h}_{i}dxdt=0 ~,\forall \hat{h}_{i}\in H_{i}.
\end{equation}
Since $y=\Lambda_{1}\bar{h}_{1}+\Lambda_{2}\bar{h}_{2}+u(g)$, and~$\tilde{y}_{d}^{i}={y}_{d}^{i}-u(g)\chi_{O_{i,d}}$, we have $\Lambda_{1}\bar{h}_{1}+\Lambda_{2}\bar{h}_{2}-\tilde{y}_{d}^{i}=y-y_{d}^{i}$ in $O_{i,d}\times(0,T)$. Note that $\hat{y}^{i}=\Lambda_{i}\hat{h}_{i}$, by (\ref{z3.13}), we get
\begin{equation}\label{z3.14}
\mu_{i}\int_{\omega_{i}\times(0,T)}\rho_{*}^{2}\bar{h}_{i}\hat{h}_{i}~dxdt+\alpha_{i}\int_{O_{i,d}\times(0,T)}(y_{1}-y_{1,d}^{i})\cdot\hat{y}_{1}^{i}+(y_{2}-y_{2,d}^{i})\cdot\hat{y}_{2}^{i}~dxdt=0,~\forall \hat{h}_{i}\in H_{i},
\end{equation}
where~$\hat{y}^{i}=(\hat{y}^{i}_{1},\hat{y}^{i}_{2})$~is the solution of~(\ref{z3.1}). By (\ref{z3.14}), we have the following expression of the Nash
equilibrium.

\begin{proposition}\label{p1}
Suppose that~$\mu_{i}~(i=1,2)$~are sufficiently large,  then for all~$g\in L^{2}(\omega\times(0,T))$, there exists a unique Nash equilibrium~$(\bar{h}_{1},\bar{h}_{2})$~for~$(J_{1},J_{2})$~such that~(\ref{1.4})~and~(\ref{1.5})~hold, simultaneously. Moreover,
\begin{eqnarray}\label{z3.15}
 \bar{h}_{i}=-\frac{1}{\mu_{i}}\rho_{*}^{-2}\varphi_{1}^{i}\chi_{\omega_{i}},~i=1,2,
\end{eqnarray}
where $(y_1, y_2, \varphi_{1}^{i},\varphi_{2}^{i})$~is the solution of the following coupled system corresponding to a function $g$:
\begin{equation}\label{z3.16}
\left\{\begin{array}{ll}
y_{1,t}+\Delta^{2}y_{1}+a_{11}y_{1}+a_{12}y_{2}=\!\!g\chi_{\omega}\!\!-\!\frac{1}{\mu_{1}}\rho_{*}^{-2}\!\varphi_{1}^{1}\chi_{\omega_{1}}\!\!-\!\frac{1}{\mu_{2}}\rho_{*}^{-2}\!\varphi_{1}^{2}\chi_{\omega_{2}}&\mbox{in}~Q,\\[3mm]
y_{2,t}+\Delta^{2}y_{2}+a_{21}y_{1}+a_{22}y_{2}=0&\mbox{in}~Q,\\[3mm]
-\varphi_{1,t}^{i}+\Delta^{2}\varphi_{1}^{i}+a_{11}\varphi_{1}^{i}+a_{21}\varphi_{2}^{i}=\alpha_{i}(y_{1}-y_{1,d}^{i})\chi_{O_{i,d}}&\mbox{in}~Q,\\[3mm]
-\varphi_{2,t}^{i}+\Delta^{2}\varphi_{2}^{i}+a_{12}\varphi_{1}^{i}+a_{22}\varphi_{2}^{i}=\alpha_{i}(y_{2}-y_{2,d}^{i})\chi_{O_{i,d}}&\mbox{in}~Q,\\[3mm]
y_{j}=\Delta y_{j}=0,~\varphi_{j}^{i}=\Delta\varphi_{j}^{i}=0&\mbox{on}~\Sigma,~j=1,2,\\[3mm]
y_{j}(x,0)=y_{j}^{0}(x),~\varphi_{j}^{i}(x,T)=0&\mbox{in}~\Omega,~j=1,2.
\end{array}\right.
\end{equation}
\end{proposition}

\textbf{Proof:}
We introduce the following adjoint system of~(\ref{z3.1}):
\begin{equation}\label{z3.17}
\left\{\begin{array}{ll}
-\varphi_{1,t}^{i}+\Delta^{2}\varphi_{1}^{i}+a_{11}\varphi_{1}^{i}+a_{21}\varphi_{2}^{i}=\alpha_{i}(y_{1}-y_{1,d}^{i})\chi_{O_{i,d}}& \mbox{ in }~Q,\\[3mm]
-\varphi_{2,t}^{i}+\Delta^{2}\varphi_{2}^{i}+a_{12}\varphi_{1}^{i}+a_{22}\varphi_{2}^{i}=\alpha_{i}(y_{2}-y_{2,d}^{i})\chi_{O_{i,d}}& \mbox{ in }~Q,\\[3mm]
\varphi_{j}^{i}=\Delta\varphi_{j}^{i}=0&\mbox{ on }~\Sigma,~j=1,2,\\[3mm]
\varphi_{j}^{i}(x,T)=0&\mbox{ in }~\Omega,~j=1,2,
\end{array}\right.
\end{equation}
where~$y=(y_{1},y_{2})$~is the solution of (\ref{1.1}),~$y_{d}^{i}=(y_{1,d}^{i},y_{2,d}^{i})$~are given functions. Multiplying the first and second equations of~(\ref{z3.17})~by~$\hat{y}_{1}^{i},\ \hat{y}_{2}^{i}$, and integrating them on $Q$, 
\begin{eqnarray*}
&\int_{Q}\varphi_{1}^{i}\cdot\hat{y}_{1,t}^{i}+\varphi_{1}^{i}\cdot\Delta^{2}\hat{y}_{1}^{i}+a_{11}\varphi_{1}^{i}\hat{y}_{1}^{i}+a_{21}\varphi_{2}^{i}\hat{y}_{1}^{i}~dxdt=\int_{Q}\alpha_{i}(y_{1}-y_{1,d}^{i})\cdot\hat{y}_{1}^{i}\chi_{O_{i,d}}~dxdt,\\[2mm]
&\int_{Q}\varphi_{2}^{i}\cdot\hat{y}_{2,t}^{i}+\varphi_{2}^{i}\cdot\Delta^{2}\hat{y}_{2}^{i}+a_{12}\varphi_{1}^{i}\hat{y}_{2}^{i}+a_{22}\varphi_{2}^{i}\hat{y}_{2}^{i}~dxdt=\int_{Q}\alpha_{i}(y_{2}-y_{2,d}^{i})\cdot\hat{y}_{2}^{i}\chi_{O_{i,d}}~dxdt.
\end{eqnarray*}
By the first and second equations of~(\ref{z3.1}), we have
\begin{eqnarray*}
&\int_{Q}\varphi_{1}^{i}(\hat{h}_{i}\chi_{\omega_{i}}-a_{12}\hat{y}_{2}^{i})+a_{21}\varphi_{2}^{i}\hat{y}_{1}^{i}~dxdt=\int_{Q}\alpha_{i}(y_{1}-y_{1,d}^{i})\cdot\hat{y}_{1}^{i}\chi_{O_{i,d}}~dxdt,\\[2mm]
&\int_{Q}\varphi_{2}^{i}(-a_{21}\hat{y}_{1}^{i})+a_{12}\varphi_{1}^{i}\hat{y}_{2}^{i}~dxdt=\int_{Q}\alpha_{i}(y_{2}-y_{2,d}^{i})\cdot\hat{y}_{2}^{i}\chi_{O_{i,d}}~dxdt.
\end{eqnarray*}
Adding up the above expressions, we have 
\begin{equation}\label{z3.18}
\int_{Q}\varphi_{1}^{i}\hat{h}_{i}\chi_{\omega_{i}}~dxdt=\alpha_{i}\int_{O_{i,d}\times(0,T)}(y_{1}-y_{1,d}^{i})\cdot\hat{y}_{1}^{i}+(y_{2}-y_{2,d}^{i})\cdot\hat{y}_{2}^{i}~dxdt.
\end{equation}
Combining this with (\ref{z3.14}), we find that 
\begin{equation*}
 \mu_{i}\int_{\omega_{i}\times(0,T)}\rho_{*}^{2}\bar{h}_{i}\hat{h}_{i}~dxdt+\int_{\omega_{i}\times(0,T)}\varphi_{1}^{i}\hat{h}_{i}~dxdt=0, ~\forall\hat{h}_{i}\in H_{i},~i=1,2.
 \end{equation*}
This implies
$$\mu_{i}\rho_{*}^{2}\bar{h}_{i}+\varphi_{1}^{i}\chi_{\omega_{i}}=0,~i=1,2,$$
which completes the proof. \endpf

\section{Reduction of the null controllability problem}

By Proposition \ref{p1}, we know that, for any control function $g$, there exists a unique Nash equilibrium pair $(\bar{h}_{1},\bar{h}_{2})$. Moreover, the solution of system~(\ref{1.1})~corresponding to the above control $g$ and $(\bar{h}_{1},\bar{h}_{2})$ is also the solution of (\ref{z3.16}). Therefore, the null controllability of the system (\ref{1.1}) can be translated into finding a control $\bar{g}\in L^{2}(\omega\times(0,T))$ such that the solution to (\ref{z3.16}) satisfies
\begin{equation}\label{z4.1}
y_{1}(\cdot,T;\bar{g})=0,~y_{2}(\cdot,T;\bar{g})=0~~\mbox{in}~\Omega,
\end{equation}
and
\begin{equation}\label{z4.2}
J(\bar{g})=\min_{g} J(g),~~\forall g\in L^{2}(\omega\times(0,T)).
\end{equation}

We consider the following adjoint equations of~(\ref{z3.16}):
\begin{equation}\label{z4.3}
\left\{\begin{array}{ll}
-\psi_{1,t}+\Delta^{2}\psi_{1}+a_{11}\psi_{1}+a_{21}\psi_{2}=\alpha_{1}\gamma_{1}^{1}\chi_{O_{1,d}}+\alpha_{2}\gamma_{1}^{2}\chi_{O_{2,d}}& \mbox{ in }~Q,\\[3mm]
-\psi_{2,t}+\Delta^{2}\psi_{2}+a_{12}\psi_{1}+a_{22}\psi_{2}=\alpha_{1}\gamma_{2}^{1}\chi_{O_{1,d}}+\alpha_{2}\gamma_{2}^{2}\chi_{O_{2,d}}& \mbox{ in }~Q,\\[3mm]
\gamma_{1,t}^{i}+\Delta^{2}\gamma_{1}^{i}+a_{11}\gamma_{1}^{i}+a_{12}\gamma_{2}^{i}=-\frac{1}{\mu_{i}}\rho_{*}^{-2}\psi_{1}\chi_{\omega_{i}}& \mbox{ in }~Q,\\[3mm]
\gamma_{2,t}^{i}+\Delta^{2}\gamma_{2}^{i}+a_{21}\gamma_{1}^{i}+a_{22}\gamma_{2}^{i}=0& \mbox{ in }~Q,\\[3mm]
\psi_{j}=\Delta \psi_{j}=0,~\gamma_{j}^{i}=\Delta\gamma_{j}^{i}=0&\mbox{ on }~\Sigma,~j=1,2,\\[3mm]
\psi_{j}(x,T)=\psi_{j}^{T},~\gamma_{j}^{i}(x,0)=0&\mbox{ in }~\Omega,~j=1,2,
\end{array}\right.
\end{equation}
where~$\psi^{T}=(\psi_{1}^{T},\psi_{2}^{T})\in L^{2}(\Omega)^{2}.$
Assume that~$O_{1,d}=O_{2,d}=O_{d}$,  then (\ref{z4.3})  can be simplified as follows
\begin{equation}\label{z4.4}
\left\{\begin{array}{ll}
-\psi_{1,t}+\Delta^{2}\psi_{1}+a_{11}\psi_{1}+a_{21}\psi_{2}=(\alpha_{1}\gamma_{1}^{1}+\alpha_{2}\gamma_{1}^{2})\chi_{O_{d}}& \mbox{ in }~Q,\\[3mm]
-\psi_{2,t}+\Delta^{2}\psi_{2}+a_{12}\psi_{1}+a_{22}\psi_{2}=(\alpha_{1}\gamma_{2}^{1}+\alpha_{2}\gamma_{2}^{2})\chi_{O_{d}}& \mbox{ in }~Q,\\[3mm]
\gamma_{1,t}^{i}+\Delta^{2}\gamma_{1}^{i}+a_{11}\gamma_{1}^{i}+a_{12}\gamma_{2}^{i}=-\frac{1}{\mu_{i}}\rho_{*}^{-2}\psi_{1}\chi_{\omega_{i}}& \mbox{ in }~Q,\\[3mm]
\gamma_{2,t}^{i}+\Delta^{2}\gamma_{2}^{i}+a_{21}\gamma_{1}^{i}+a_{22}\gamma_{2}^{i}=0& \mbox{ in }~Q,\\[3mm]
\psi_{j}=\Delta \psi_{j}=0,~\gamma_{j}^{i}=\Delta\gamma_{j}^{i}=0&\mbox{ on }~\Sigma,~j=1,2,\\[3mm]
\psi_{j}(x,T)=\psi_{j}^{T},~\gamma_{j}^{i}(x,0)=0&\mbox{ in }~\Omega,~j=1,2.
\end{array}\right.
\end{equation}

By the duality theory, we can transform the null controllability of the system (\ref{z3.16}) with respect to $y$ into the following observability inequality.

\begin{proposition}\label{p8}
If there exists a  constant~$C>0$ such that the solutions $\psi$ and $\gamma^{i}=(\gamma_{1}^{i},\gamma_{2}^{i})~(i=1,2)$~of~(\ref{z4.4}) satisfy
\begin{eqnarray}\label{z4.5}
\begin{array}{rl}
&\displaystyle\int_{\Omega}|\psi_{1}(x,0)|^{2}~dx+\int_{\Omega}|\psi_{2}(x,0)|^{2}~dx+\sum_{i=1}^{2}\int_{Q}|\gamma_{1}^{i}|^{2}~dxdt+\sum_{i=1}^{2}\int_{Q}|\gamma_{2}^{i}|^{2}~dxdt
\\[3mm]&\leq \tilde{C}\displaystyle\int_{\omega\times(0,T)}|\psi_{1}|^{2}~dxdt,
\ \forall\ \psi^{T}\in L^{2}(\Omega)^{2}, 
\end{array}
\end{eqnarray}
then, the solutions of (\ref{z3.16}) satisfy (\ref{z4.1}) and (\ref{z4.2}).
\end{proposition}

\textbf{Proof.} For any $\psi^{T}\in L^{2}(\Omega)^{2}$, we construct the following  functional
\begin{eqnarray*}
\begin{array}{rl}
&\tilde{F}(\psi^{T})=\displaystyle\frac{1}{2}\int_{\omega\times(0,T)}|\psi_{1}|^{2}~dxdt+\int_{\Omega}y_{1}^{0}\psi_{1}(x,0)~dx+\int_{\Omega}y_{2}^{0}\psi_{2}(x,0)~dx
\\[3mm]&\quad\quad\quad\quad-\displaystyle\sum_{i=1}^{2}\int_{O_{d}\times(0,T)}\alpha_{i}\gamma_{1}^{i}y_{1,d}^{i}~dxdt-\sum_{i=1}^{2}\int_{O_{d}\times(0,T)}\alpha_{i}\gamma_{2}^{i}y_{2,d}^{i}~dxdt.
\end{array}
\end{eqnarray*}
It is easy to check that~$\tilde{F}$~is continuous and strictly
convex. Next, we prove that $\tilde{F}$~is coercive.

By Cauchy inequality with~$\varepsilon$, we have
\begin{eqnarray*}
\begin{array}{rl}
&\!\!\!\!\!\tilde{F}(\psi^{T})\geq\displaystyle\frac{1}{2}\int_{\omega\times(0,T)}|\psi_{1}|^{2}~dxdt
\\[3mm]&\quad\quad\quad\!-\displaystyle\frac{\varepsilon}{2}\left[\!\int_{\Omega}\!|\psi_{1}(x,0)|^{2}dx\!\!+\!\!\!\int_{\Omega}\!|\psi_{2}(x,0)|^{2}dx\!\!+\!\!\!\sum_{i=1}^{2}\!\!\int_{O_{d}\times(0,T)}\!\!\!|\gamma_{1}^{i}|^{2}dxdt\!\!+\!\!\!\sum_{i=1}^{2}\!\!\int_{O_{d}\times(0,T)}\!\!\!|\gamma_{2}^{i}|^{2}dxdt\right]
\\[3mm]&\quad\quad\quad\!-\displaystyle\frac{1}{2\varepsilon}\left[\!\int_{\Omega}|y_{1}^{0}|^{2}dx\!\!+\!\!\!\int_{\Omega}|y_{2}^{0}|^{2}dx\!\!+\!\!\!\sum_{i=1}^{2}\!\!\int_{O_{d}\times(0,T)}\!\!\alpha_{i}^{2}|y_{1,d}^{i}|^{2}dxdt\!\!+\!\!\!\sum_{i=1}^{2}\!\!\int_{O_{d}\times(0,T)}\!\!\alpha_{i}^{2}|y_{2,d}^{i}|^{2}dxdt\right].
\end{array}
\end{eqnarray*}
By (\ref{z4.5}), we deduce that
\begin{eqnarray*}
\begin{array}{rl}
&\!\!\!\!\!\tilde{F}(\psi^{T})\displaystyle\geq\left(\frac{1}{2}-\frac{\varepsilon \tilde{C}}{2}\right)\int_{\omega\times(0,T)}|\psi_{1}|^{2}~dxdt
\\[3mm]&\quad\quad\quad\!-\displaystyle\frac{1}{2\varepsilon}\left[\!\int_{\Omega}|y_{1}^{0}|^{2}dx\!\!+\!\!\!\int_{\Omega}|y_{2}^{0}|^{2}dx\!\!+\!\!\!\sum_{i=1}^{2}\!\!\int_{O_{d}\times(0,T)}\!\!\alpha_{i}^{2}|y_{1,d}^{i}|^{2}dxdt\!\!+\!\!\!\sum_{i=1}^{2}\!\!\int_{O_{d}\times(0,T)}\!\!\alpha_{i}^{2}|y_{2,d}^{i}|^{2}dxdt\right].
\end{array}
\end{eqnarray*}
Take~$\varepsilon=1/2\tilde{C}$, then
\begin{eqnarray*}
\begin{array}{rl}
&\!\!\!\!\!\tilde{F}(\psi^{T})\displaystyle\geq\frac{1}{4}\int_{\omega\times(0,T)}|\psi_{1}|^{2}~dxdt
\\[3mm]&\quad\quad\quad\!-\tilde{C}\displaystyle\left[\!\int_{\Omega}|y_{1}^{0}|^{2}dx\!\!+\!\!\!\int_{\Omega}|y_{2}^{0}|^{2}dx\!\!+\!\!\!\sum_{i=1}^{2}\!\!\int_{O_{d}\times(0,T)}\!\!\alpha_{i}^{2}|y_{1,d}^{i}|^{2}dxdt\!\!+\!\!\!\sum_{i=1}^{2}\!\!\int_{O_{d}\times(0,T)}\!\!\alpha_{i}^{2}|y_{2,d}^{i}|^{2}dxdt\right].
\end{array}
\end{eqnarray*}
Combining (\ref{1.8}) with the above inequality we get that $\tilde{F}$~is coercive. Therefore, there exists a minimizer~$\hat{\psi}^{T}$~such that~$\tilde{F}(\hat{\psi}^{T})=\displaystyle\min_{\psi^{T}\in L^{2}(\Omega)^{2}}\tilde{F}(\psi^{T})$. 
Denote by  $\hat{\psi}=(\hat{\psi}_{1},\hat{\psi}_{2})$~and~$\hat{\gamma}^{i}=(\hat{\gamma}^{i}_{1},\hat{\gamma}^{i}_{2})$ the solutions of (\ref{z4.4}) corresponding to $\hat{\psi}^{T}$, and denote by $\psi$~and~$\gamma^{i}=(\gamma_{1}^{i},\gamma_{2}^{i})$ the solutions of (\ref{z4.4}) corresponding to  ~$\psi^{T}\in L^{2}(\Omega)^{2}$. Since (\ref{z4.4}) is a linear system, it follows that 
\begin{eqnarray*}
\begin{array}{rl}
&\!\!\!\!\!\displaystyle\lim_{\beta\rightarrow0}\frac{\tilde{F}(\hat{\psi}^{T}+\beta\psi^{T})-\tilde{F}(\hat{\psi}^{T})}{\beta}
=0.
\end{array}
\end{eqnarray*}
By calculation, this implies that 
\begin{eqnarray}\label{tt}
\begin{array}{rl}
&\displaystyle\int_{\omega\times(0,T)}\hat{\psi}_{1}\psi_{1}~dxdt-\displaystyle\sum_{i=1}^{2}\int_{O_{d}\times(0,T)}\alpha_{i}\gamma_{1}^{i}y_{1,d}^{i}~dxdt-\sum_{i=1}^{2}\int_{O_{d}\times(0,T)}\alpha_{i}\gamma_{2}^{i}y_{2,d}^{i}~dxdt
\\[3mm]&\quad\quad\quad+\displaystyle\int_{\Omega}y_{1}^{0}\psi_{1}(x,0)~dx+\displaystyle\int_{\Omega}y_{2}^{0}\psi_{2}(x,0)~dx=0.
\end{array}
\end{eqnarray}

On the other hand, multiplying both sides of the first four equations of (\ref{z3.16}) by $\psi_1,\psi_2,\gamma_1^i$ and $\gamma_2^i$, respectively, and integrating them on $Q$, by (\ref{z3.16}) and (\ref{z4.3}), we obtain that
\begin{eqnarray}\label{g5.5}
\begin{array}{rl}
&\displaystyle\int_{Q}\left[(\alpha_{1}\gamma_{1}^{1}\chi_{O_{1,d}}+\alpha_{2}\gamma_{1}^{2}\chi_{O_{2,d}})\cdot y_{1}+(\alpha_{1}\gamma_{2}^{1}\chi_{O_{1,d}}+\alpha_{2}\gamma_{2}^{2}\chi_{O_{2,d}})\cdot y_{2}\right]dxdt
\\[2mm]&\quad+\displaystyle\int_{\Omega}y_{1}\psi_{1} dx\bigg{|}^{T}_{0}+\int_{\Omega}y_{2}\psi_{2} dx\bigg{|}^{T}_{0}
\\[2mm]&=\displaystyle\int_{Q}(g\chi_{\omega}-\frac{1}{\mu_{1}}\rho_{*}^{-2}\varphi_{1}^{1}\chi_{\omega_{1}}-\frac{1}{\mu_{2}}\rho_{*}^{-2}\varphi_{1}^{2}\chi_{\omega_{2}})\cdot\psi_{1}dxdt,
\end{array}
\end{eqnarray}
\begin{eqnarray}\label{g5.6}
\begin{array}{rl}
\displaystyle\int_{Q}\left(-\frac{1}{\mu_{i}}\rho_{*}^{-2}\varphi_{1}^{i}\psi_{1}\chi_{\omega_{i}}\right)dxdt=\!\displaystyle\int_{Q}\left[\alpha_{i}(y_{1}-y_{1,d}^{i})\gamma_{1}^{i}\chi_{O_{i,d}}\!+\!\alpha_{i}(y_{2}-y_{2,d}^{i})\gamma_{2}^{i}\chi_{O_{i,d}}\right]dxdt.
\end{array}
\end{eqnarray}
By (\ref{g5.5}) and (\ref{g5.6}), we have
\begin{eqnarray*}
\begin{array}{rl}
&\quad\displaystyle\int_{\omega\times(0,T)}g\psi_{1} ~dxdt-\displaystyle\sum_{i=1}^{2}\displaystyle\int_{O_{i,d}\times(0,T)}\alpha_{i}\gamma_{1}^{i}y_{1,d}^{i}~dxdt-\sum_{i=1}^{2}\displaystyle\int_{O_{i,d}\times(0,T)}\alpha_{i}\gamma_{2}^{i}y_{2,d}^{i}~dxdt
\\[2mm]
&\quad=\displaystyle\int_{\Omega}y_{1}\psi_{1}~dx\bigg{|}^{T}_{0}+\displaystyle\int_{\Omega}y_{2}\psi_{2}~dx\bigg{|}^{T}_{0}.
\end{array}
\end{eqnarray*}
This implies that, if (\ref{z4.1}) holds, then
\begin{eqnarray}\label{g5.7}
\begin{array}{rl}
&\displaystyle\int_{\omega\times(0,T)}g\psi_{1} ~dxdt-\sum_{i=1}^{2}\int_{O_{i,d}\times(0,T)}\alpha_{i}\gamma_{1}^{i}y_{1,d}^{i}~dxdt-\sum_{i=1}^{2}\int_{O_{i,d}\times(0,T)}\alpha_{i}\gamma_{2}^{i}y_{2,d}^{i}~dxdt
\\[3mm]&\quad+\displaystyle\int_{\Omega}y_{1}^{0}\psi_{1}(x,0)~dx+\displaystyle\int_{\Omega}y_{2}^{0}\psi_{2}(x,0)~dx=0.
\end{array}
\end{eqnarray}
Combining (\ref{tt}) with (\ref{g5.7}),  we know that the control can be taken as $\bar{g}=\hat{\psi}_{1}\chi_{\omega}$, making (\ref{z4.1}) and (\ref{z4.2}) hold.
\endpf

\section{Carleman estimate of the coupled fourth order parabolic system}

In order to establish the Carleman estimate of the coupled fourth order parabolic system, we need to introduce some weight functions:
\begin{equation}\label{z2.3}
\sigma(x,t)=\displaystyle\frac{e^{4\lambda\|\eta\|_{L^{\infty}(\Omega)}}-e^{\lambda\left(2\|\eta\|_{L^{\infty}(\Omega)}+\eta(x)\right)}}{t^{1/2}(T-t)^{1/2}},
\end{equation}
\begin{equation}\label{z2.4} 
\tau(x,t)=\frac{e^{\lambda\left(2\|\eta\|_{L^{\infty}(\Omega)}+\eta(x)\right)}}{t^{1/2}(T-t)^{1/2}}.
\end{equation}
Here~$\eta\in C^{4}(\overline{\Omega})$ satisfies
\begin{equation}\label{z2.5}
\eta|_{\partial\Omega}=0,~|\nabla\eta|\geq C_{0}>0~in ~\Omega\setminus\overline{\omega^{\prime}},~\omega^{\prime}\subset\omega.
\end{equation}
Set
\begin{equation}\label{z3.6}
\sigma^{*}(t)=\max_{\overline{\Omega}}\sigma(x,t),~\displaystyle\rho_{*}(t)=\displaystyle e^{\frac{s\sigma^{*}}{2}}.
\end{equation}
We recall a known global Carleman estimate for a single fourth-order parabolic equation.

\begin{lemma}\label{l1}\upcite{2019guerrero}
There exist constants~$C>0$ and $s_0,\lambda_0\geq 1$, such that for any~$\lambda\geq \lambda_{0}$,~$s\geq s_{0}(T^{1/2}+T)$ and any solution $\varphi$ of 
\begin{equation}\label{z2.2}
\left\{\begin{array}{ll}
-\varphi_{t}+\Delta^{2}\varphi=f& \mbox{ in }~Q,\\[3mm]
\varphi=\Delta\varphi=0&\mbox{ on }~\Sigma,\\[3mm]
\varphi(x,T)=\varphi_{T}(x)&\mbox{ in }~\Omega,
\end{array}\right.
\end{equation}
it holds that 
\begin{eqnarray}\label{z2.6}
\begin{array}{rl}
&\displaystyle\int_{Q}e^{-2s\sigma}\Big{[}\lambda^{8}(s\tau)^{6}|\varphi|^{2}+\lambda^{6}(s\tau)^{4}|\nabla\varphi|^{2}+\lambda^{4}(s\tau)^{3}|\Delta\varphi|^{2}+\lambda^{4}(s\tau)^{2}|\nabla^{2}\varphi|^{2}
\\[3mm]&\quad+\lambda^{2}(s\tau)|\nabla\Delta\varphi|^{2}+s^{-1}\tau^{-1}(|\varphi_{t}|^{2}+|\Delta^{2}\varphi|^{2})\Big{]}dxdt
\\[3mm]&\leq C_{0}\left(\lambda^{8}\displaystyle\int_{\omega\times(0,T)}e^{-2s\sigma}(s\tau)^{7}|\varphi|^{2}dxdt+\int_{Q}e^{-2s\sigma}|f|^{2}dxdt\right),
\end{array}
\end{eqnarray}
where~$\varphi_{T}\in L^{2}(\Omega), f\in L^{2}(Q).$
\end{lemma}

Let~$\theta_{j}=\alpha_{1}\gamma_{j}^{1}+\alpha_{2}\gamma_{j}^{2}~(j=1,2)$. Then (\ref{z4.4}) can be simplified as follows:
\begin{equation}\label{z4.6}
\left\{\begin{array}{ll}
-\psi_{1,t}+\Delta^{2}\psi_{1}+a_{11}\psi_{1}+a_{21}\psi_{2}=\theta_{1}\chi_{O_{d}}& \mbox{ in }~Q,\\[3mm]
-\psi_{2,t}+\Delta^{2}\psi_{2}+a_{12}\psi_{1}+a_{22}\psi_{2}=\theta_{2}\chi_{O_{d}}& \mbox{ in }~Q,\\[3mm]
\theta_{1,t}+\Delta^{2}\theta_{1}+a_{11}\theta_{1}+a_{12}\theta_{2}=-\rho_{*}^{-2}(\frac{\alpha_{1}}{\mu_{1}}\chi_{\omega_{1}}+\frac{\alpha_{2}}{\mu_{2}}\chi_{\omega_{2}})\psi_{1}& \mbox{ in }~Q,\\[3mm]
\theta_{2,t}+\Delta^{2}\theta_{2}+a_{21}\theta_{1}+a_{22}\theta_{2}=0& \mbox{ in }~Q,\\[3mm]
\psi_{j}=\Delta \psi_{j}=0,~\theta_{j}=\Delta\theta_{j}=0&\mbox{ on }~\Sigma,~j=1,2,\\[3mm]
\psi_{j}(x,T)=\psi_{j}^{T},~\theta_{j}(x,0)=0&\mbox{ in }~\Omega,~j=1,2,
\end{array}\right.
\end{equation}
where~$\psi_{j}^{T}\in L^{2}(\Omega),~j=1,2.$~

In what follows, we put
\begin{eqnarray*}
\begin{array}{rl}
&I(\varphi):=\displaystyle\int_{Q}e^{-2s\sigma}\Big{[}\lambda^{8}(s\tau)^{6}|\varphi|^{2}+\lambda^{6}(s\tau)^{4}|\nabla\varphi|^{2}+\lambda^{4}(s\tau)^{3}|\Delta\varphi|^{2}+\lambda^{4}(s\tau)^{2}|\nabla^{2}\varphi|^{2}
\\[3mm]&\quad\quad\quad\quad+\lambda^{2}(s\tau)|\nabla\Delta\varphi|^{2}+s^{-1}\tau^{-1}(|\varphi_{t}|^{2}+|\Delta^{2}\varphi|^{2})\Big{]}dxdt,
\end{array}
\end{eqnarray*}
and
\begin{equation}\label{z4.7}
\rho_{0}=\displaystyle\min_{t\in[0,T]}\rho_{*}(t), \ \text{where} \ \rho_{*}(t)\ \text{is given in} \ (\ref{z3.6}).
\end{equation}

We derive the following Carleman estimate for~(\ref{z4.6}) in this section.

\begin{proposition}\label{p3}
Suppose that~$O_{d}\cap \omega\neq\emptyset$,~$\mu_{i}$~are sufficiently large and $a_{21}\geq a_{0}>0$~or~$-a_{21}\geq a_{0}>0$ in $(O_{d}\cap\omega)\times(0,T)$. Then there exist constants~$C_{1}>0$~and~$\lambda_{1},s_{1}\geq1$, such that for any~$\lambda\geq\lambda_{1},s\geq s_{1}$, the solution of system (\ref{z4.6}) satisfies
\begin{eqnarray}\label{z4.8}
\begin{array}{rl}
I(\psi_{1})+I(\psi_{2})+I(\theta_{1})+I(\theta_{2})\leq C_{1}\lambda^{24}\displaystyle\int_{\omega\times(0,T)}e^{-2s\sigma}(s\tau)^{34}|\psi_{1}|^{2}dxdt,
\end{array}
\end{eqnarray}
where~$C_{1}=C_{1}(\Omega,\omega,\rho_{0},\alpha_{1},\alpha_{2},\mu_{1},\mu_{2}).$~
\end{proposition}

{\bf Proof.} Since~$O_{d}\cap\omega\neq\emptyset$, we choose a nonempty open set~$\omega^{\prime}$~satisfying~$\overline{\omega^{\prime}}\subset O_{d}\cap\omega.$~

First, applying Lemma~\ref{l1}~to the solutions~$\psi_{j}$~and~$\theta_{j}$~of (\ref{z4.6}), respectively, we obtain
\begin{equation*}
 I(\psi_{1})\leq C_{0}\left(\lambda^{8}\displaystyle\int_{\omega^{\prime}\times(0,T)}e^{-2s\sigma}(s\tau)^{7}|\psi_{1}|^{2}dxdt+\int_{Q}e^{-2s\sigma}|\theta_{1}\chi_{O_{d}}-a_{11}\psi_{1}-a_{21}\psi_{2}|^{2}dxdt\right),
 \end{equation*}
\begin{equation*}
 I(\psi_{2})\leq C_{0}\left(\lambda^{8}\displaystyle\int_{\omega^{\prime}\times(0,T)}e^{-2s\sigma}(s\tau)^{7}|\psi_{2}|^{2}dxdt+\int_{Q}e^{-2s\sigma}|\theta_{2}\chi_{O_{d}}-a_{12}\psi_{1}-a_{22}\psi_{2}|^{2}dxdt\right),
 \end{equation*}
\begin{equation*}
I(\theta_{1})\leq C_{0}\left(\lambda^{8}\!\!\!\!\displaystyle\int_{\omega^{\prime}\times(0,T)}\!\!\!\!\!\!e^{-2s\sigma}(s\tau)^{7}|\theta_{1}|^{2}\!dxdt\!\!+\!\!\!\!\displaystyle\int_{Q}\!\!e^{-2s\sigma}\Big{|}\!\!-\!\!\rho_{*}^{-2}(\frac{\alpha_{1}}{\mu_{1}}\chi_{\omega_{1}}\!\!+\!\!\frac{\alpha_{2}}{\mu_{2}}\chi_{\omega_{2}})\psi_{1}\!\!-\!a_{11}\theta_{1}\!\!-\!a_{12}\theta_{2}\Big{|}^{2}\!\!dxdt\!\right),
\end{equation*}
and
\begin{equation*}
 I(\theta_{2})\leq C_{0}\left(\lambda^{8}\displaystyle\int_{\omega^{\prime}\times(0,T)}e^{-2s\sigma}(s\tau)^{7}|\theta_{2}|^{2}dxdt+\int_{Q}e^{-2s\sigma}\Big{|}-a_{21}\theta_{1}-a_{22}\theta_{2}\Big{|}^{2}dxdt\right).
\end{equation*}
Taking a sufficiently large $s$, it follows that 
\begin{eqnarray}\label{z4.9}
\begin{array}{rl}
&I(\psi_{1})+I(\psi_{2})+I(\theta_{1})+I(\theta_{2})
\\[3mm]&\leq C_{0}\lambda^{8}\displaystyle\int_{\omega^{\prime}\times(0,T)}e^{-2s\sigma}(s\tau)^{7}(|\psi_{1}|^{2}+|\psi_{2}|^{2}+|\theta_{1}|^{2}+|\theta_{2}|^{2})dxdt.
\end{array}
\end{eqnarray}

Next, we estimate the last two terms in the right-hand side of  (\ref{z4.9}).
Choose a nonempty set~$\omega^{\prime\prime}$ satisfying
~$\overline{\omega^{\prime}}\subset\omega^{\prime\prime}\subset O_{d}\cap\omega,$ and take a cut-off function~$\xi\in C^{\infty}(\mathbb{R}^{N})$ such that 
\\$0\leq\xi\leq1$~in~$\mathbb{R}^{N}$,
$\xi=1$~in~$\omega^{\prime}$,~$supp\xi\subseteq\omega^{\prime\prime}$, and
$$\frac{\nabla\xi}{\xi^{1/2}}\in L^{\infty}(\Omega)^{N},~\frac{\Delta\xi}{\xi^{1/2}}\in L^{\infty}(\Omega),~\frac{\nabla^{2}\xi}{\xi^{1/2}}\in L^{\infty}(\Omega)^{N^{2}},~\frac{\nabla\Delta\xi}{\xi^{1/2}}\in L^{\infty}(\Omega)^{N},~\frac{\Delta^{2}\xi}{\xi^{1/2}}\in L^{\infty}(\Omega).$$
We define~$\delta=e^{-2s\sigma}(s\tau)^{7}$. By~(\ref{z4.6})~and 
$\xi|_{\partial\Omega}=0,~\nabla\xi|_{\partial\Omega}=0,~\delta|_{t=0}=0,~\delta|_{t=T}=0$, we get 
\begin{eqnarray}\label{z4.10}
\begin{array}{rl}
&C_{0}\lambda^{8}\displaystyle\int_{\omega^{\prime}\times(0,T)}e^{-2s\sigma}(s\tau)^{7}|\theta_{1}|^{2}dxdt+C_{0}\lambda^{8}\displaystyle\int_{\omega^{\prime}\times(0,T)}e^{-2s\sigma}(s\tau)^{7}|\theta_{2}|^{2}dxdt
\\[3mm]&\leq C_{0}\lambda^{8}\displaystyle\int_{O_{d}\times(0,T)}\delta\xi|\theta_{1}|^{2}dxdt+C_{0}\lambda^{8}\displaystyle\int_{O_{d}\times(0,T)}\delta\xi|\theta_{2}|^{2}dxdt
\\[3mm]&= C_{0}\lambda^{8}\displaystyle\int_{O_{d}\times(0,T)}\delta\xi\theta_{1}(-\psi_{1,t}+\Delta^{2}\psi_{1}+a_{11}\psi_{1}+a_{21}\psi_{2})dxdt
\\[3mm]&\quad+C_{0}\lambda^{8}\displaystyle\int_{O_{d}\times(0,T)}\delta\xi\theta_{2}(-\psi_{2,t}+\Delta^{2}\psi_{2}+a_{12}\psi_{1}+a_{22}\psi_{2})dxdt
\\[3mm]&=
C_{0}\lambda^{8}\displaystyle\int_{Q}\left[\psi_{1}(\delta\xi\theta_{1})_{t}+\psi_{1}\Delta^{2}(\delta\xi\theta_{1})+\psi_{1}\delta\xi a_{11}\theta_{1}+\psi_{2}\delta\xi a_{21}\theta_{1}\right]dxdt
\\[3mm]&\quad+C_{0}\lambda^{8}\displaystyle\int_{Q}\left[\psi_{2}(\delta\xi\theta_{2})_{t}+\psi_{2}\Delta^{2}(\delta\xi\theta_{2})+\psi_{1}\delta\xi a_{12}\theta_{2}+\psi_{2}\delta\xi a_{22}\theta_{2}\right]dxdt
\\[3mm]&=C_{0}\lambda^{8}\!\!\displaystyle\int_{Q}\!\psi_{1}\delta\xi(\theta_{1,t}\!+\!\Delta^{2}\theta_{1}\!+\!a_{11}\theta_{1}\!+\!a_{12}\theta_{2})\!+\!\psi_{2}\delta\xi(\theta_{2,t}\!+\!\Delta^{2}\theta_{2}\!+\!a_{21}\theta_{1}\!+\!a_{22}\theta_{2})dxdt
\\[3mm]&\quad+C_{0}\lambda^{8}\!\displaystyle\sum_{j=1}^{2}\!\!\int_{Q}\psi_{j}\big{[}(\delta\xi)_{t}\theta_{j}\!+\!\Delta^{2}(\delta\xi)\theta_{j}\!+\!4\nabla\Delta(\delta\xi)\!\cdot\!\nabla\theta_{j}\!+\!2\Delta(\delta\xi)\!\cdot\!\Delta\theta_{j}\!+\!2\nabla^{2}(\delta\xi)\!\cdot\!\Delta\theta_{j}
\\[3mm]&\quad\quad\quad\quad\quad\quad\quad+2\Delta(\delta\xi)\cdot\nabla^{2}\theta_{j}+4\nabla(\delta\xi)\cdot\nabla\Delta\theta_{j}\big{]}dxdt.
\end{array}
\end{eqnarray}

On the one hand,  it is calculated by (\ref{z2.3}) and (\ref{z2.4}) that
\begin{eqnarray}\label{z4.11}
\begin{array}{rl}
&|\delta_{t}|\leq C e^{-2s\sigma}(s\tau)^{10},\quad|\nabla\delta|\leq C\lambda e^{-2s\sigma}(s\tau)^{8},\quad
|\Delta\delta|\leq C\lambda^{2}e^{-2s\sigma}(s\tau)^{9},\quad
\\[2mm]&|\nabla^{2}\delta|\leq C\lambda^{2} e^{-2s\sigma}(s\tau)^{9},
|\nabla\Delta\delta|\leq C\lambda^{3}e^{-2s\sigma}(s\tau)^{10},\quad|\Delta^{2}\delta|\leq C\lambda^{4} e^{-2s\sigma}(s\tau)^{11}.
\end{array}
\end{eqnarray}
On the other hand, we obtain from (\ref{z4.6}) that
\begin{eqnarray*}
\begin{array}{rl}
&I_{1}:=\displaystyle\int_{Q}\psi_{1}\delta\xi(\theta_{1,t}+\Delta^{2}\theta_{1}+a_{11}\theta_{1}+a_{12}\theta_{2})~dxdt
\\[3mm]&\quad=\displaystyle\int_{Q}\psi_{1}\delta\xi\rho_{*}^{-2}(-\frac{\alpha_{1}}{\mu_{1}}\chi_{\omega_{1}}-\frac{\alpha_{2}}{\mu_{2}}\chi_{\omega_{2}})\psi_{1}~dxdt
\\[3mm]&\quad\displaystyle\leq C_{1}\int_{\omega^{\prime\prime}\times(0,T)}\psi_{1}^{2}\delta ~dxdt = C_{1}\int_{\omega^{\prime\prime}\times(0,T)}e^{-2s\sigma}(s\tau)^{7}|\psi_{1}|^{2} ~dxdt,
\end{array}
\end{eqnarray*}
and
\begin{eqnarray*}
\begin{array}{rl}
&I_{2}:=\displaystyle\int_{Q}\psi_{2}\delta\xi(\theta_{2,t}+\Delta^{2}\theta_{2}+a_{21}\theta_{1}+a_{22}\theta_{2})~dxdt=0.
\end{array}
\end{eqnarray*}
Further, by (\ref{z4.11}) and H\"{o}lder's inequality, for any $\varepsilon>0$, we have
\begin{eqnarray*}
\begin{array}{rl}
&I_{3}:=\displaystyle\int_{Q}\psi_{j}(\delta\xi)_{t}\theta_{j}~dxdt=\int_{Q}\psi_{j}\theta_{j}\xi\delta_{t}~dxdt
\\[3mm]&\quad\displaystyle\leq \int_{\omega^{\prime\prime}\times(0,T)}\psi_{j}\theta_{j}\delta_{t} ~dxdt
\leq C\int_{\omega^{\prime\prime}\times(0,T)}\psi_{j}\theta_{j} e^{-2s\sigma}(s\tau)^{10} ~dxdt
\\[3mm]&\quad \displaystyle\leq \frac{C}{2\varepsilon}\int_{\omega^{\prime\prime}\times(0,T)}e^{-2s\sigma}(s\tau)^{14}|\psi_{j}|^{2}dxdt+\frac{C\varepsilon}{2}\int_{\omega^{\prime\prime}\times(0,T)}e^{-2s\sigma}(s\tau)^{6}|\theta_{j}|^{2}dxdt,
\end{array}
\end{eqnarray*}
\begin{eqnarray*}
\begin{array}{rl}
&I_{4}:=\displaystyle\int_{Q}\psi_{j}\theta_{j}\Delta^{2}(\delta\xi)~dxdt
\\[3mm]&\quad=\displaystyle\int_{Q}\psi_{j}\theta_{j}(\Delta^{2}\delta\cdot\xi+4\nabla\Delta\delta\cdot\nabla\xi+2\Delta\delta\cdot\Delta\xi+2\nabla^{2}\delta\cdot\Delta\xi+2\Delta\delta\cdot\nabla^{2}\xi
\\[3mm]&\quad\quad\quad\quad\quad+4\nabla\delta\cdot\nabla\Delta\xi+\delta\cdot\Delta^{2}\xi)~dxdt
\\[3mm]&\quad=\displaystyle\int_{Q}\psi_{j}\theta_{j}(\Delta^{2}\delta\cdot\xi+4\nabla\Delta\delta\cdot\frac{\nabla\xi}{\xi^{\frac{1}{2}}}\cdot\xi^{\frac{1}{2}}+2\Delta\delta\cdot\frac{\Delta\xi}{\xi^{\frac{1}{2}}}\cdot\xi^{\frac{1}{2}}+2\nabla^{2}\delta\cdot\frac{\Delta\xi}{\xi^{\frac{1}{2}}}\cdot\xi^{\frac{1}{2}}
\\[3mm]&\quad\quad\quad\quad\quad  +2\displaystyle\Delta\delta\cdot\frac{\nabla^{2}\xi}{\xi^{\frac{1}{2}}}\cdot\xi^{\frac{1}{2}}+4\displaystyle\nabla\delta\cdot\frac{\nabla\Delta\xi}{\xi^{\frac{1}{2}}}\cdot\xi^{\frac{1}{2}}+\delta\cdot\frac{\Delta^{2}\xi}{\xi^{\frac{1}{2}}}\cdot\xi^{\frac{1}{2}})~dxdt
\\[3mm]&\quad\leq C\displaystyle\int_{\omega^{\prime\prime}\times(0,T)}\psi_{j}\theta_{j}(|\Delta^{2}\delta|+|\nabla\Delta\delta|+|\Delta\delta|+|\nabla^{2}\delta|+|\nabla\delta|+|\delta|)dxdt
\\[3mm]&\quad\leq \displaystyle\frac{\lambda^{8}C}{2\varepsilon}\int_{\omega^{\prime\prime}\times(0,T)}e^{-2s\sigma}(s\tau)^{16}|\psi_{j}|^{2}dxdt+\frac{C\varepsilon}{2}\int_{\omega^{\prime\prime}\times(0,T)}e^{-2s\sigma}(s\tau)^{6}|\theta_{j}|^{2}dxdt,
\end{array}
\end{eqnarray*}
\begin{eqnarray*}
\begin{array}{rl}
&I_{5}:=\displaystyle\int_{Q}\psi_{j}\cdot\nabla\theta_{j}\cdot\nabla\Delta(\delta\xi)~dxdt
\\[3mm]&\quad\leq\displaystyle\frac{\lambda^{6}C}{2\varepsilon}\int_{\omega^{\prime\prime}\times(0,T)}e^{-2s\sigma}(s\tau)^{16}|\psi_{j}|^{2}dxdt+\frac{C\varepsilon}{2}\int_{\omega^{\prime\prime}\times(0,T)}e^{-2s\sigma}(s\tau)^{4}|\nabla\theta_{j}|^{2}dxdt,
\end{array}
\end{eqnarray*}
\begin{eqnarray*}
\begin{array}{rl}
&I_{6}:=\displaystyle\int_{Q}\psi_{j}\cdot\Delta\theta_{j}\cdot\Delta(\delta\xi)~dxdt
\\[3mm]&\quad\leq\displaystyle\frac{\lambda^{4}C}{2\varepsilon}\int_{\omega^{\prime\prime}\times(0,T)}e^{-2s\sigma}(s\tau)^{15}|\psi_{j}|^{2}dxdt+\frac{C\varepsilon}{2}\int_{\omega^{\prime\prime}\times(0,T)}e^{-2s\sigma}(s\tau)^{3}|\Delta\theta_{j}|^{2}dxdt,
\end{array}
\end{eqnarray*}
\begin{eqnarray*}
\begin{array}{rl}
&I_{7}:=\displaystyle\int_{Q}\psi_{j}\cdot\Delta\theta_{j}\cdot\nabla^{2}(\delta\xi)~dxdt
\\[3mm]&\quad\leq\displaystyle\frac{\lambda^{4}C}{2\varepsilon}\int_{\omega^{\prime\prime}\times(0,T)}e^{-2s\sigma}(s\tau)^{15}|\psi_{j}|^{2}dxdt+\frac{C\varepsilon}{2}\int_{\omega^{\prime\prime}\times(0,T)}e^{-2s\sigma}(s\tau)^{3}|\Delta\theta_{j}|^{2}dxdt,
\end{array}
\end{eqnarray*}
\begin{eqnarray*}
\begin{array}{rl}
&I_{8}:=\displaystyle\int_{Q}\psi_{j}\cdot\nabla^{2}\theta_{j}\cdot\Delta(\delta\xi)~dxdt
\\[3mm]&\quad\leq\displaystyle\frac{\lambda^{4}C}{2\varepsilon}\int_{\omega^{\prime\prime}\times(0,T)}e^{-2s\sigma}(s\tau)^{16}|\psi_{j}|^{2}dxdt+\frac{C\varepsilon}{2}\int_{\omega^{\prime\prime}\times(0,T)}e^{-2s\sigma}(s\tau)^{2}|\nabla^{2}\theta_{j}|^{2}dxdt,
\end{array}
\end{eqnarray*}
and
\begin{eqnarray*}
\begin{array}{rl}
&I_{9}:=\displaystyle\int_{Q}\psi_{j}\cdot\nabla\Delta\theta_{j}\cdot\nabla(\delta\xi)~dxdt
\\[3mm]&\quad\leq\displaystyle\frac{\lambda^{2}C}{2\varepsilon}\int_{\omega^{\prime\prime}\times(0,T)}e^{-2s\sigma}(s\tau)^{15}|\psi_{j}|^{2}dxdt+\frac{C\varepsilon}{2}\int_{\omega^{\prime\prime}\times(0,T)}e^{-2s\sigma}(s\tau)|\nabla\Delta\theta_{j}|^{2}dxdt.
\end{array}
\end{eqnarray*}
Substituting the above inequalities of $I_{1}$ to $I_{9}$ into~(\ref{z4.10}), by (\ref{z4.9}), and taking $\varepsilon$
sufficiently small, we obtain
\begin{eqnarray}\label{z4.13}
\begin{array}{rl}
&I(\psi_{1})+I(\psi_{2})+I(\theta_{1})+I(\theta_{2})
\\[3mm]&\leq C_{1}\lambda^{16} \displaystyle\int_{\omega^{\prime\prime}\times(0,T)}e^{-2s\sigma}(s\tau)^{16}|\psi_{1}|^{2}dxdt
+C_{0}\lambda^{16}\displaystyle\int_{\omega^{\prime\prime}\times(0,T)}e^{-2s\sigma}(s\tau)^{16}|\psi_{2}|^{2}dxdt.
\end{array}
\end{eqnarray}

Now, we deal with the last term in (\ref{z4.13}).
Let $\tilde{\omega}$ be a subset of $O_{d}\cap\omega$, such that~$\overline{\omega^{\prime\prime}}\subset\tilde{\omega}\subset O_{d}\cap\omega$. 
Take a cut-off function~$\tilde{\xi}\in C^{\infty}(\mathbb{R}^{N})$~satisfying $0\leq\tilde{\xi}\leq1$, in~$\mathbb{R}^{N}$, $\tilde{\xi}=1$, in $\omega^{\prime\prime}$,~$supp\tilde{\xi}\subseteq\tilde{\omega}$, and
$$\frac{\nabla\tilde{\xi}}{\tilde{\xi}^{1/2}}\in L^{\infty}(\Omega)^{N},~\frac{\Delta\tilde{\xi}}{\tilde{\xi}^{1/2}}\in L^{\infty}(\Omega),~\frac{\nabla^{2}\tilde{\xi}}{\tilde{\xi}^{1/2}}\in L^{\infty}(\Omega)^{N^{2}},~\frac{\nabla\Delta\tilde{\xi}}{\tilde{\xi}^{1/2}}\in L^{\infty}(\Omega)^{N},~\frac{\Delta^{2}\tilde{\xi}}{\tilde{\xi}^{1/2}}\in L^{\infty}(\Omega).$$
Define $\tilde{\delta}=e^{-2s\sigma}(s\tau)^{16}$. Notice that  $a_{21}\geq a_{0}>0$ in $(O_{d}\cap\omega)\times(0,T)$ and $\tilde{\xi}|_{\partial\Omega}=0,~\nabla\tilde{\xi}|_{\partial\Omega}=0,~\tilde{\delta}|_{t=0}=0,~\tilde{\delta}|_{t=T}=0$. By (\ref{z4.6}),  with  the same argument as in (\ref{z4.10}), we have
\begin{eqnarray}\label{z4.14}
\begin{array}{rl}
&a_{0}\displaystyle\int_{\omega^{\prime\prime}\times(0,T)}e^{-2s\sigma}(s\tau)^{16}|\psi_{2}|^{2}~dxdt=a_{0}\displaystyle\int_{\omega^{\prime\prime}\times(0,T)}\tilde{\delta}\tilde{\xi}|\psi_{2}|^{2}~dxdt
\\[3mm]&\leq\displaystyle\int_{\omega^{\prime\prime}\times(0,T)}\tilde{\delta}\tilde{\xi}a_{21}|\psi_{2}|^{2}~dxdt\leq\displaystyle\int_{Q}\tilde{\delta}\tilde{\xi}\psi_{2}(\psi_{1,t}-\Delta^{2}\psi_{1}-a_{11}\psi_{1}+\theta_{1}\chi_{O_{d}})~dxdt
\\[3mm]&=\displaystyle\int_{Q}\left[-\psi_{1}(\tilde{\delta}\tilde{\xi}\psi_{2})_{t}-\psi_{1}\Delta^{2}(\tilde{\delta}\tilde{\xi}\psi_{2})-a_{11}\tilde{\delta}\tilde{\xi}\psi_{1}\psi_{2}+\tilde{\delta}\tilde{\xi}\theta_{1}\psi_{2}\chi_{O_{d}}\right]~dxdt.
\end{array}
\end{eqnarray}
Note that 
\begin{eqnarray}\label{z4.15}
\begin{array}{rl}
&|\tilde{\delta}_{t}|\leq C e^{-2s\sigma}(s\tau)^{19},\quad|\nabla\tilde{\delta}|\leq C\lambda e^{-2s\sigma}(s\tau)^{17},
\quad|\Delta\tilde{\delta}|\leq C\lambda^{2} e^{-2s\sigma}(s\tau)^{18},
\\[3mm]
&|\nabla^{2}\tilde{\delta}|\leq C\lambda^{2} e^{-2s\sigma}(s\tau)^{18},
\quad|\nabla\Delta\tilde{\delta}|\leq C\lambda^{3} e^{-2s\sigma}(s\tau)^{19},\quad|\Delta^{2}\tilde{\delta}|\leq C\lambda^{4} e^{-2s\sigma}(s\tau)^{20}.
\end{array}
\end{eqnarray}
Denote by $K_{1}$~to~$K_{4}$ the items in the right hand side of~(\ref{z4.14}), respectively. By the definition of $\tilde{\xi}$,
 H\"{o}lder's inequality and (\ref{z4.15}), we have
\begin{eqnarray*}
\begin{array}{rl}
&K_{1}\leq\displaystyle\int_{Q}\big{|}\psi_{1}(\tilde{\delta}\tilde{\xi}\psi_{2})_{t}\big{|}~dxdt=\displaystyle\int_{Q}\big{|}\psi_{1}\tilde{\xi}(\tilde{\delta}_{t}\psi_{2}+\tilde{\delta}\psi_{2,t})\big{|}~dxdt
\\[3mm]&\quad\leq\displaystyle\int_{\tilde{\omega}\times(0,T)}\big{|}\psi_{1}\psi_{2}\tilde{\delta}_{t}+\psi_{1}\psi_{2,t}\tilde{\delta}\big{|}~dxdt
\\[3mm]&\quad\leq\displaystyle\frac{C}{2\varepsilon}\int_{\tilde{\omega}\times(0,T)}e^{-2s\sigma}(s\tau)^{33}|\psi_{1}|^{2}~dxdt+\displaystyle\frac{C\varepsilon}{2}\int_{\tilde{\omega}\times(0,T)}e^{-2s\sigma}(s\tau)^{6}|\psi_{2}|^{2}~dxdt
\\[3mm]&\quad\quad+\displaystyle\frac{C\varepsilon}{2}\int_{\tilde{\omega}\times(0,T)}e^{-2s\sigma}(s\tau)^{-1}|\psi_{2,t}|^{2}~dxdt,
\end{array}
\end{eqnarray*}
\begin{eqnarray*}
&&K_{2}\leq\displaystyle\int_{Q}\big{|}\psi_{1}\Delta^{2}(\tilde{\delta}\tilde{\xi}\psi_{2})\big{|}~dxdt
\\[3mm]&&\quad=\displaystyle\int_{Q}\psi_{1}(\Delta^{2}\big{(}\tilde{\delta}\tilde{\xi})\cdot\psi_{2}+4\nabla\Delta(\tilde{\delta}\tilde{\xi})\cdot\nabla\psi_{2}+2\Delta(\tilde{\delta}\tilde{\xi})\cdot\Delta\psi_{2}+2\nabla^{2}(\tilde{\delta}\tilde{\xi})\cdot\Delta\psi_{2}
\\[3mm]&&\quad\quad\quad\quad\quad+2\Delta(\tilde{\delta}\tilde{\xi})\cdot\nabla^{2}\psi_{2}+4\nabla(\tilde{\delta}\tilde{\xi})\cdot\nabla\Delta\psi_{2}+\tilde{\delta}\tilde{\xi}\cdot\Delta^{2}\psi_{2}\big{)}~dxdt
\\[3mm]&&\quad\leq\displaystyle\frac{\lambda^{8}C}{2\varepsilon}\int_{\tilde{\omega}\times(0,T)}e^{-2s\sigma}(s\tau)^{34}|\psi_{1}|^{2}~dxdt+\displaystyle\frac{C\varepsilon}{2}\int_{\tilde{\omega}\times(0,T)}e^{-2s\sigma}(s\tau)^{6}|\psi_{2}|^{2}~dxdt
\\[3mm]&&\quad\quad+\displaystyle\frac{C\varepsilon}{2}\int_{\tilde{\omega}\times(0,T)}e^{-2s\sigma}(s\tau)^{4}|\nabla\psi_{2}|^{2}~dxdt+\displaystyle\frac{C\varepsilon}{2}\int_{\tilde{\omega}\times(0,T)}e^{-2s\sigma}(s\tau)^{3}|\Delta\psi_{2}|^{2}~dxdt
\\[3mm]&&\quad\quad+\displaystyle\frac{C\varepsilon}{2}\int_{\tilde{\omega}\times(0,T)}e^{-2s\sigma}(s\tau)^{2}|\nabla^{2}\psi_{2}|^{2}~dxdt+\displaystyle\frac{C\varepsilon}{2}\int_{\tilde{\omega}\times(0,T)}e^{-2s\sigma}(s\tau)|\nabla\Delta\psi_{2}|^{2}~dxdt
\\[3mm]&&\quad\quad+\displaystyle\frac{C\varepsilon}{2}\int_{\tilde{\omega}\times(0,T)}e^{-2s\sigma}(s\tau)^{-1}|\Delta^{2}\psi_{2}|^{2}~dxdt,
\end{eqnarray*}
\begin{eqnarray*}
\begin{array}{rl}
&K_{3}\leq\displaystyle\int_{Q}\big{|}a_{11}\tilde{\delta}\tilde{\xi}\psi_{1}\psi_{2}\big{|}~dxdt
\\[3mm]&\quad\leq\displaystyle\frac{C}{2\varepsilon}\int_{\tilde{\omega}\times(0,T)}e^{-2s\sigma}(s\tau)^{26}|\psi_{1}|^{2}~dxdt+\displaystyle\frac{C\varepsilon}{2}\int_{\tilde{\omega}\times(0,T)}e^{-2s\sigma}(s\tau)^{6}|\psi_{2}|^{2}~dxdt,
\end{array}
\end{eqnarray*}
and
\begin{eqnarray*}
\begin{array}{rl}
&K_{4}\leq\displaystyle\int_{Q}\big{|}\tilde{\delta}\tilde{\xi}\theta_{1}\psi_{2}\chi_{O_{d}}\big{|}~dxdt\leq\displaystyle\int_{\tilde{\omega}\times(0,T)}\big{|}e^{-2s\sigma}(s\tau)^{16}\theta_{1}\psi_{2}\big{|}~dxdt
\\[3mm]&\quad\leq\displaystyle\frac{1}{2}\int_{\tilde{\omega}\times(0,T)}e^{-4s\sigma}(s\tau)^{32}|\psi_{2}|^{2}~dxdt+\displaystyle\frac{1}{2}\int_{\tilde{\omega}\times(0,T)}|\theta_{1}|^{2}~dxdt.
\end{array}
\end{eqnarray*}
We take a  sufficiently large $s$, such that $e^{-2s\sigma}(s\tau)^{26}<\tilde{\varepsilon}<1$, then 
$$K_{4}\leq\frac{\tilde{\varepsilon}}{2}\int_{\tilde{\omega}\times(0,T)}e^{-2s\sigma}(s\tau)^{6}|\psi_{2}|^{2}~dxdt+\frac{1}{2}\int_{\tilde{\omega}\times(0,T)}|\theta_{1}|^{2}~dxdt.$$
Similarly, the above estimates still hold when~$-a_{21}\geq a_{0}>0$ in $(O_{d}\cap\omega)\times(0,T)$.
 We substitute the estimates of $K_{1},K_{2},K_{3}$ and $ K_{4}$~into~(\ref{z4.14}). By (\ref{z4.13}), we get
\begin{eqnarray}\label{z4.16}
\begin{array}{rl}
&I(\psi_{1})+I(\psi_{2})+I(\theta_{1})+I(\theta_{2})
\\[3mm]&\leq C_{1}\lambda^{24}\displaystyle\int_{\tilde{\omega}\times(0,T)}e^{-2s\sigma}(s\tau)^{34}|\psi_{1}|^{2}dxdt+C_{0}\lambda^{16}\displaystyle\int_{\tilde{\omega}\times(0,T)}|\theta_{1}|^{2}dxdt.
\end{array}
\end{eqnarray}
We will finally estimate the local term concerning $\theta_{1}$. By the classical energy estimates for the third and fourth equations in~(\ref{z4.6}), we have
\begin{equation}\label{z4.17}
\displaystyle\int_{Q}(|\theta_{1}|^{2}+|\theta_{2}|^{2})~dxdt\leq C_{2}(\frac{\alpha_{1}^{2}}{\mu_{1}^{2}}+\frac{\alpha_{2}^{2}}{\mu_{2}^{2}})\displaystyle\int_{Q}|\rho_{*}^{-2}\psi_{1}|^{2}~dxdt,
\end{equation}
where
$C_{2}=C_{2}(\Omega,\omega,T,\|a_{11}\|_{L^{\infty}(Q)},\|a_{12}\|_{L^{\infty}(Q)},\|a_{21}\|_{L^{\infty}(Q)},\|a_{22}\|_{L^{\infty}(Q)}).$~

Since~$\rho_{*}^{-4}\leq e^{-2s\sigma}$, take $\mu_{1},\mu_{2}$ sufficiently large, then the right term of (\ref{z4.17}) can be absorbed by~$I(\psi_{1})$, by (\ref{z4.16}) and (\ref{z4.17}), it holds that 
\begin{eqnarray}\label{z4.18}
\begin{array}{rl}
I(\psi_{1})+I(\psi_{2})+I(\theta_{1})+I(\theta_{2})\leq C_{1}\lambda^{24}\displaystyle\int_{\tilde{\omega}\times(0,T)}e^{-2s\sigma}(s\tau)^{34}|\psi_{1}|^{2}dxdt,
\end{array}
\end{eqnarray}
then the desired estimate in  Proposition~\ref{p3} follows.
\endpf

\section{Observability of the coupled fourth order parabolic system}

\begin{proposition}\label{p2}
Suppose that~$O_{d}\cap \omega\neq\emptyset$,~$\mu_{i}$~are sufficiently large and $a_{21}\geq a_{0}>0$~or~$-a_{21}\geq a_{0}>0$ in $(O_{d}\cap\omega)\times(0,T)$. Then there exists a positive constant~$\tilde{C}$, such that for all~$\psi^{T}\in L^{2}(\Omega)^{2}$, the solutions $\psi$ and $\gamma^{i}=(\gamma_{1}^{i},\gamma_{2}^{i})~(i=1,2)$ of (\ref{z4.4}) satisfy
\begin{eqnarray}\label{z4.5}
\begin{array}{rl}
&\displaystyle\int_{\Omega}|\psi_{1}(x,0)|^{2}~dx+\int_{\Omega}|\psi_{2}(x,0)|^{2}~dx+\sum_{i=1}^{2}\int_{Q}|\gamma_{1}^{i}|^{2}~dxdt+\sum_{i=1}^{2}\int_{Q}|\gamma_{2}^{i}|^{2}~dxdt
\\[3mm]&\leq \tilde{C}\displaystyle\int_{\omega\times(0,T)}|\psi_{1}|^{2}~dxdt.
\end{array}
\end{eqnarray}
\end{proposition}

\textbf{Proof.} Let $\tilde{\eta}\in C^{1}[0,T]$ be a function satisfying:
$$0\leq\tilde{\eta}\leq1\ \text{in}\ [0,T],\ \tilde{\eta}=1\ \text{in}\ [0,T/2],\ \tilde{\eta}=0 \ \text{in}\ [3T/4,T],\ |\tilde{\eta}_{t}|\leq C/T.$$
Multiplying both sides of the first two equations of  (\ref{z4.6}) by $\tilde{\eta}\psi_{1}$ and $\tilde{\eta}\psi_{2}$, respectively, and integrating them in~$\Omega$, by H\"{o}lder's inequality, we have
\begin{eqnarray*}
\begin{array}{rl}
&\displaystyle-\frac{1}{2}\frac{d}{dt}\int_{\Omega}\tilde{\eta}\cdot(\psi_{1}^{2}+\psi_{2}^{2})dx+\displaystyle\int_{\Omega}\tilde{\eta}\cdot(|\Delta\psi_{1}|^{2}+|\Delta\psi_{2}|^{2})dx
\\[3mm]&=-\displaystyle\int_{\Omega}\tilde{\eta} a_{11}|\psi_{1}|^{2}dx-\displaystyle\int_{\Omega}\tilde{\eta} a_{22}|\psi_{2}|^{2}dx-\displaystyle\int_{\Omega}\tilde{\eta}(a_{12}+a_{21})\psi_{1}\psi_{2}dx
\\[3mm]&\quad-\displaystyle\frac{1}{2}\int_{\Omega}\tilde{\eta}_{t}(\psi_{1}^{2}+\psi_{2}^{2})dx+\int_{\Omega}\theta_{1}\tilde{\eta}\psi_{1}\chi_{O_{d}}dx+\int_{\Omega}\theta_{2}\tilde{\eta}\psi_{2}\chi_{O_{d}}dx
\\[3mm]&\leq\displaystyle\|a_{11}\|_{L^{\infty}(Q)}\int_{\Omega}\tilde{\eta} |\psi_{1}|^{2}dx+\displaystyle\|a_{22}\|_{L^{\infty}(Q)}\int_{\Omega}\tilde{\eta}|\psi_{2}|^{2}dx
\\[3mm]&\quad+\displaystyle\frac{\|a_{12}\|_{L^{\infty}(Q)}+\|a_{21}\|_{L^{\infty}(Q)}}{2}\left(\displaystyle\int_{\Omega}\tilde{\eta}|\psi_{1}|^{2}dx+\int_{\Omega}\tilde{\eta}|\psi_{2}|^{2}dx\right)-\displaystyle\frac{1}{2}\int_{\Omega}\tilde{\eta}_{t}(\psi_{1}^{2}+\psi_{2}^{2})dx
\\[3mm]&\quad+\displaystyle\frac{1}{2}\int_{\Omega}\tilde{\eta}(\psi_{1}^{2}+\psi_{2}^{2})dx+\displaystyle\frac{1}{2}\int_{\Omega}\tilde{\eta}(\theta_{1}^{2}+\theta_{2}^{2})dx.
\end{array}
\end{eqnarray*}
This implies that 
\begin{eqnarray}\label{z4.19}
\begin{array}{rl}
&-\displaystyle\frac{d}{dt}\int_{\Omega}\tilde{\eta}\cdot(\psi_{1}^{2}+\psi_{2}^{2})dx+2\displaystyle\int_{\Omega}\tilde{\eta}\cdot(|\Delta\psi_{1}|^{2}+|\Delta\psi_{2}|^{2})dx
\\[3mm]&\leq(1+2A)\displaystyle\int_{\Omega}\tilde{\eta}(\psi_{1}^{2}+\psi_{2}^{2})dx+\displaystyle\int_{\Omega}\tilde{\eta}(\theta_{1}^{2}+\theta_{2}^{2})dx-\displaystyle\int_{\Omega}\tilde{\eta}_{t}(\psi_{1}^{2}+\psi_{2}^{2})dx,
\end{array}
\end{eqnarray}
where~$A:=\displaystyle\sum_{i=1}^{2}\displaystyle\sum_{j=1}^{2}\|a_{ij}\|_{L^{\infty}(Q)}.$ Multiplying both sides of (\ref{z4.19})~by~$e^{(1+2A)t}$~and integrating it in~$[0,T]$, we obtain
\begin{eqnarray}\label{z4.20}
\begin{array}{rl}
&\displaystyle\int_{0}^{T}-\frac{d}{dt}\left(\displaystyle e^{(1+2A)t}\!\int_{\Omega}\tilde{\eta}(\psi_{1}^{2}+\psi_{2}^{2})dx\right)dt+2\int_{0}^{T}\displaystyle\int_{\Omega}\tilde{\eta} e^{(1+2A)t}(|\Delta\psi_{1}|^{2}+|\Delta\psi_{2}|^{2})dxdt
\\[3mm]&\leq \displaystyle\int_{0}^{T}\int_{\Omega}\tilde{\eta} e^{(1+2A)t}(\theta_{1}^{2}+\theta_{2}^{2}) dxdt-\displaystyle\int_{0}^{T}\int_{\Omega}\tilde{\eta}_{t}e^{(1+2A)t} (\psi_{1}^{2}+\psi_{2}^{2})dxdt.
\end{array}
\end{eqnarray}
Noticing that 
$$|\tilde{\eta}_{t}|\leq C/T~in~[T/2,3T/4];~\tilde{\eta}_{t}=0~in~[0,T/2)\cup(3T/4,T],$$
and by (\ref{z4.20}),~we get
\begin{eqnarray*}
\begin{array}{rl}
&-e^{(1+2A)t}\displaystyle\int_{\Omega}\tilde{\eta}(\psi_{1}^{2}+\psi_{2}^{2})dx\bigg{|}^{T}_{0}+2\displaystyle\int_{0}^{\frac{T}{2}}\int_{\Omega}e^{(1+2A)t}(|\Delta\psi_{1}|^{2}+|\Delta\psi_{2}|^{2})dxdt
\\[3mm]&\quad+2\displaystyle\int_{\frac{T}{2}}^{\frac{3}{4}T}\int_{\Omega}\tilde{\eta} e^{(1+2A)t}(|\Delta\psi_{1}|^{2}+|\Delta\psi_{2}|^{2})dxdt
\\[3mm]&\leq\displaystyle\int_{0}^{\frac{3}{4}T}\int_{\Omega}e^{(1+2A)t}(\theta_{1}^{2}+\theta_{2}^{2})dxdt+\frac{C}{T}\displaystyle\int_{\frac{T}{2}}^{\frac{3}{4}T}\int_{\Omega}e^{(1+2A)t}(\psi_{1}^{2}+\psi_{2}^{2})dxdt.
\end{array}
\end{eqnarray*}
Since~$e^{(1+2A)t}$~have lower and upper bounds in $[0,T]$, we have
\begin{eqnarray}\label{z4.21}
\begin{array}{rl}
&\displaystyle\int_{\Omega}(|\psi_{1}(x,0)|^{2}+|\psi_{2}(x,0)|^{2})dx+\displaystyle\int_{0}^{\frac{T}{2}}\int_{\Omega}(|\Delta\psi_{1}|^{2}+|\Delta\psi_{2}|^{2})dxdt
\\[3mm]&\leq C_{2}\left(\displaystyle\int_{0}^{\frac{3}{4}T}\int_{\Omega}(\theta_{1}^{2}+\theta_{2}^{2})dxdt+ \displaystyle\int_{\frac{T}{2}}^{\frac{3}{4}T}\int_{\Omega}(\psi_{1}^{2}+\psi_{2}^{2})dxdt\right).
\end{array}
\end{eqnarray}

On the other hand, for any $\tau\in [0,T/2]$, we multiply both sides of (\ref{z4.19})~by~$e^{(1+2A)t}$~and integrate it in~$[\tau,T]$, then 
\begin{eqnarray}\label{z4.22}
\begin{array}{rl}
&\displaystyle \int_{\tau}^{T}-\frac{d}{dt}\left(\displaystyle e^{(1+2A)t}\!\int_{\Omega}\tilde{\eta}(\psi_{1}^{2}+\psi_{2}^{2})dx\right)dt+2\int_{\tau}^{T}\displaystyle\int_{\Omega}\tilde{\eta} e^{(1+2A)t}(|\Delta\psi_{1}|^{2}+|\Delta\psi_{2}|^{2})dxdt
\\[3mm]&\leq \displaystyle\int_{\tau}^{T}\int_{\Omega}\tilde{\eta} e^{(1+2A)t}(\theta_{1}^{2}+\theta_{2}^{2}) dxdt-\displaystyle\int_{\tau}^{T}\int_{\Omega}\tilde{\eta}_{t}e^{(1+2A)t}(\psi_{1}^{2}+\psi_{2}^{2})dxdt.
\end{array}
\end{eqnarray}
Since~$\tilde{\eta}(T)=0$,~$\tilde{\eta}(\tau)=1$, and~$e^{(1+2A)t}$~is bounded in~$[0,T]$, we obtain
\begin{eqnarray}\label{z4.23}
\begin{array}{rl}
&\displaystyle\int_{\Omega}(|\psi_{1}(x,\tau)|^{2}+|\psi_{2}(x,\tau)|^{2})dx
\\[3mm]&\leq C_{2}\left(\displaystyle\int_{0}^{\frac{3}{4}T}\int_{\Omega}(|\theta_{1}|^{2}+|\theta_{2}|^{2})dxdt+ \int_{\frac{T}{2}}^{\frac{3}{4}T}\displaystyle\int_{\Omega}(|\psi_{1}|^{2}+|\psi_{2}|^{2})dxdt\right).
\end{array}
\end{eqnarray}
Integrating~(\ref{z4.23})~with respect to~$\tau$~in~$[0,T/2]$~, we arrive at 
\begin{eqnarray}\label{z4.24}
\begin{array}{rl}
&\displaystyle\int_{0}^{\frac{T}{2}}\int_{\Omega}(|\psi_{1}|^{2}+|\psi_{2}|^{2})dxdt
\\[3mm]&\leq C_{2}\left(\displaystyle\int_{0}^{\frac{3}{4}T}\int_{\Omega}(|\theta_{1}|^{2}+|\theta_{2}|^{2})dxdt+ \int_{\frac{T}{2}}^{\frac{3}{4}T}\displaystyle\int_{\Omega}(|\psi_{1}|^{2}+|\psi_{2}|^{2})dxdt\right).
\end{array}
\end{eqnarray}
Combining (\ref{z4.21})~with~(\ref{z4.24}), it holds that 
\begin{eqnarray}\label{z4.25}
\begin{array}{rl}
&\displaystyle\int_{\Omega}(|\psi_{1}(x,0)|^{2}+|\psi_{2}(x,0)|^{2})dx+\displaystyle\int_{0}^{\frac{T}{2}}\int_{\Omega}(|\psi_{1}|^{2}+|\psi_{2}|^{2}+|\Delta\psi_{1}|^{2}+|\Delta\psi_{2}|^{2})dxdt
\\[3mm]&\leq \displaystyle C_{2}\left(\int_{0}^{\frac{3}{4}T}\int_{\Omega}(|\theta_{1}|^{2}+|\theta_{2}|^{2})dxdt+ \displaystyle\int_{\frac{T}{2}}^{\frac{3}{4}T}\int_{\Omega}(|\psi_{1}|^{2}+|\psi_{2}|^{2})dxdt\right).
\end{array}
\end{eqnarray}
Define the function
\begin{eqnarray}\label{z4.26}
l(t)=\left\{\begin{array}{ll}
\frac{T}{2},&0\leq t\leq\frac{T}{2},\\[2mm]
\displaystyle t^{\frac{1}{2}}(T-t)^{\frac{1}{2}},&\frac{T}{2}\leq t\leq T,
\end{array}\right.
\end{eqnarray}
and the following associated weight functions
\begin{equation}\label{z4.27}
\bar{\sigma}(x,t)=\displaystyle\frac{e^{4\lambda\|\eta\|_{L^{\infty}(\Omega)}}-e^{\lambda\left(2\|\eta\|_{L^{\infty}(\Omega)}+\eta(x)\right)}}{l(t)},
\end{equation}
\begin{equation}\label{z4.28}
\quad\quad\bar{\tau}(x,t)=\frac{e^{\lambda\left(2\|\eta\|_{L^{\infty}(\Omega)}+\eta(x)\right)}}{l(t)}, \ \text{where}\ \eta\ \text{is\ as\ in}\ (\ref{z2.5}).
\end{equation}

Since~$\bar{\sigma}(x,t)$ and $\bar{\tau}(x,t)$~have lower and upper bounds in~$\overline{\Omega}\times[0,T/2]$, similar to  (\ref{z4.25}), we have
\begin{eqnarray}\label{z4.29}
\begin{array}{rl}
&\displaystyle\int_{\Omega}(|\psi_{1}(x,0)|^{2}+|\psi_{2}(x,0)|^{2})dx+\displaystyle \bar{I}_{[0,\frac{T}{2}]}(\psi_{1})+\displaystyle \bar{I}_{[0,\frac{T}{2}]}(\psi_{2})
\\[3mm]&\leq \displaystyle C_{2}\left(\int_{0}^{\frac{3}{4}T}\int_{\Omega}(|\theta_{1}|^{2}+|\theta_{2}|^{2})dxdt+ \displaystyle\int_{\frac{T}{2}}^{\frac{3}{4}T}\int_{\Omega}(|\psi_{1}|^{2}+|\psi_{2}|^{2})dxdt\right),
\end{array}
\end{eqnarray}
where~$\displaystyle\bar{I}_{[a,b]}(\psi)=\displaystyle\int_{a}^{b}\int_{\Omega}e^{-2s\bar{\sigma}}(s\bar{\tau})^{6}|\psi|^{2}dxdt+\int_{a}^{b}\int_{\Omega}e^{-2s\bar{\sigma}}(s\bar{\tau})^{3}|\Delta\psi|^{2}dxdt.$~

Adding the terms $\bar{I}_{[0,\frac{T}{2}]}(\theta_{1})$ and $\bar{I}_{[0,\frac{T}{2}]}(\theta_{2})$ to both sides of~(\ref{z4.29}), and noticing that $e^{-2s\bar{\sigma}}(s\bar{\tau})^{6}$~has a lower bound in~$\overline{\Omega}\times[0,T/2]$,  it concludes that
\begin{eqnarray}\label{z4.30}
\begin{array}{rl}
&\displaystyle\int_{\Omega}(|\psi_{1}(x,0)|^{2}+|\psi_{2}(x,0)|^{2})dx+\displaystyle \bar{I}_{[0,\frac{T}{2}]}(\psi_{1})+\displaystyle \bar{I}_{[0,\frac{T}{2}]}(\psi_{2})+\displaystyle \bar{I}_{[0,\frac{T}{2}]}(\theta_{1})+\displaystyle \bar{I}_{[0,\frac{T}{2}]}(\theta_{2})
\\[3mm]
&\leq \displaystyle C_{2}\left(\displaystyle\int_{\frac{T}{2}}^{\frac{3}{4}T}\int_{\Omega}(|\psi_{1}|^{2}+|\psi_{2}|^{2}+|\theta_{1}|^{2}+|\theta_{2}|^{2})dxdt+\displaystyle \bar{I}_{[0,\frac{T}{2}]}(\theta_{1})+\bar{I}_{[0,\frac{T}{2}]}(\theta_{2}) \right).
\end{array}
\end{eqnarray}
In order to deal with the terms~$\displaystyle \bar{I}_{[0,\frac{T}{2}]}(\theta_{1})$ and $\displaystyle \bar{I}_{[0,\frac{T}{2}]}(\theta_{2})$~in the right side of (\ref{z4.30}), we use classical energy estimate for~$\theta_{j}$~in~(\ref{z4.6}), then
\begin{equation}\label{z4.31}
\displaystyle\int_{0}^{\frac{T}{2}}\int_{\Omega}(|\theta_{1}|^{2}+|\theta_{2}|^{2}+|\Delta\theta_{1}|^{2}+|\Delta\theta_{2}|^{2})dxdt\leq C_{2}(\frac{\alpha_{1}^{2}}{\mu_{1}^{2}}+\frac{\alpha_{2}^{2}}{\mu_{2}^{2}})\displaystyle\int_{0}^{\frac{T}{2}}\int_{\Omega}|\rho_{*}^{-2}\psi_{1}|^{2}dxdt.
\end{equation}
Since~$e^{-2s\bar{\sigma}}(s\bar{\tau})^{6}$~has a upper bound in~$\overline{\Omega}\times[0,T/2]$, we get
$$\displaystyle \bar{I}_{[0,\frac{T}{2}]}(\theta_{1})+\bar{I}_{[0,\frac{T}{2}]}(\theta_{2})\leq C\displaystyle\int_{0}^{\frac{T}{2}}\int_{\Omega}(|\theta_{1}|^{2}+|\theta_{2}|^{2}+|\Delta\theta_{1}|^{2}+|\Delta\theta_{2}|^{2})dxdt.$$
By (\ref{z4.31}), (\ref{z4.7}) and the fact that~$e^{-2s\bar{\sigma}}(s\bar{\tau})^{6}$~has a lower bound in~$\overline{\Omega}\times[0,T/2]$, we arrive at 
\begin{eqnarray}\label{z4.32}
\begin{array}{rl}
&\displaystyle\bar{I}_{[0,\frac{T}{2}]}(\theta_{1})+\bar{I}_{[0,\frac{T}{2}]}(\theta_{2})\leq C_{2}\displaystyle(\frac{\alpha_{1}^{2}}{\mu_{1}^{2}}+\frac{\alpha_{2}^{2}}{\mu_{2}^{2}})\displaystyle\int_{0}^{\frac{T}{2}}\int_{\Omega}|\rho_{*}^{-2}\psi_{1}|^{2}dxdt
\\[3mm]&\quad\quad\quad\quad\quad\quad\quad\quad\quad\leq C_{2}\rho_{0}^{-4}\displaystyle(\frac{\alpha_{1}^{2}}{\mu_{1}^{2}}+\frac{\alpha_{2}^{2}}{\mu_{2}^{2}})\bar{I}_{[0,\frac{T}{2}]}(\psi_{1}).
\end{array}
\end{eqnarray}
Taking $\mu_{1},\mu_{2}$ large enough, by~(\ref{z4.30}),~(\ref{z4.32}) and Proposition \ref{p3}, we have 
\begin{eqnarray}\label{z4.33}
\begin{array}{rl}
&\displaystyle\int_{\Omega}(|\psi_{1}(x,0)|^{2}+|\psi_{2}(x,0)|^{2})dx+\displaystyle \bar{I}_{[0,\frac{T}{2}]}(\psi_{1})+\displaystyle \bar{I}_{[0,\frac{T}{2}]}(\psi_{2})+\displaystyle \bar{I}_{[0,\frac{T}{2}]}(\theta_{1})+\displaystyle \bar{I}_{[0,\frac{T}{2}]}(\theta_{2})
\\[3mm]
&\leq \displaystyle C_{2}\displaystyle\int_{\frac{T}{2}}^{\frac{3}{4}T}\int_{\Omega}(|\psi_{1}|^{2}+|\psi_{2}|^{2}+|\theta_{1}|^{2}+|\theta_{2}|^{2})dxdt+C_{2}\rho_{0}^{-4}(\frac{\alpha_{1}^{2}}{\mu_{1}^{2}}+\frac{\alpha_{2}^{2}}{\mu_{2}^{2}})\bar{I}_{[0,\frac{T}{2}]}(\psi_{1})
\\[3mm]
&\leq C_{2}\Big{(}I(\psi_{1})+I(\psi_{2})+I(\theta_{1})+I(\theta_{2})\Big{)}
\\[3mm]
&\leq\tilde{C}\lambda^{24}\displaystyle\int_{\omega\times(0,T)}e^{-2s\sigma}(s\tau)^{34}|\psi_{1}|^{2}dxdt
\leq\tilde{C}\displaystyle\int_{\omega\times(0,T)}|\psi_{1}|^{2}dxdt,
\end{array}
\end{eqnarray}
where~$\tilde{C}=\tilde{C}(\Omega,\omega,T,\|a_{11}\|_{L^{\infty}(Q)},\|a_{12}\|_{L^{\infty}(Q)},\|a_{21}\|_{L^{\infty}(Q)},\|a_{22}\|_{L^{\infty}(Q)},\rho_{0},\alpha_{1},\alpha_{2},\mu_{1},\mu_{2}).$~
\medskip

In addition, we see that~$\bar{\sigma}=\sigma$,~$\bar{\tau}=\tau$ in~$[T/2,T]$. By Proposition \ref{p3}, we deduce that
\begin{eqnarray}\label{z4.34}
\begin{array}{rl}
&\bar{I}_{[\frac{T}{2},T]}(\psi_{1})+\bar{I}_{[\frac{T}{2},T]}(\psi_{2}) +\bar{I}_{[\frac{T}{2},T]}(\theta_{1})+\bar{I}_{[\frac{T}{2},T]}(\theta_{2})
\\[3mm]&\leq I(\psi_{1})+I(\psi_{2})+I(\theta_{1})+I(\theta_{2})
\leq C_{1}\displaystyle\int_{\omega\times(0,T)}|\psi_{1}|^{2}dxdt.
\end{array}
\end{eqnarray}
Combining~(\ref{z4.33})~with~(\ref{z4.34}), we have
\begin{eqnarray}\label{z4.35}
\begin{array}{rl}
&\displaystyle\int_{\Omega}(|\psi_{1}(x,0)|^{2}+|\psi_{2}(x,0)|^{2})dx+\displaystyle \bar{I}_{[0,T]}(\psi_{1})+\displaystyle \bar{I}_{[0,T]}(\psi_{2})+\displaystyle \bar{I}_{[0,T]}(\theta_{1})+\displaystyle \bar{I}_{[0,T]}(\theta_{2})
\\[3mm]&\leq \tilde{C}\displaystyle\int_{\omega\times(0,T)}|\psi_{1}|^{2}dxdt.
\end{array}
\end{eqnarray}

On the other hand, by the classical energy estimates for~$\gamma^{i}~(i=1,2)$~in (\ref{z4.4}), we get
\begin{equation}\label{z4.36}
\displaystyle\sum_{i=1}^{2}\int_{Q}(|\gamma_{1}^{i}|^{2}+|\gamma_{2}^{i}|^{2})dxdt\leq C_{2}\displaystyle\sum_{i=1}^{2}\frac{1}{\mu_{i}^{2}}\int_{\omega_{i}\times(0,T)}|\rho_{*}^{-2}\psi_{1}|^{2}dxdt.
\end{equation}
By (\ref{z3.6}),~(\ref{z4.35}) and~$\bar{\sigma}\leq\sigma\leq\sigma^{*}$, we obtain
\begin{eqnarray}\label{z4.37}
\begin{array}{rl}
&\displaystyle\int_{\omega_{i}\times(0,T)}|\rho_{*}^{-2}\psi_{1}|^{2}dxdt\leq\displaystyle\int_{\omega_{i}\times(0,T)}e^{-2s\bar{\sigma}}|\psi_{1}|^{2}(s\bar{\tau})^{6}\frac{1}{(s\bar{\tau})^{6}}dxdt
\\[3mm]&\leq \displaystyle\max_{\overline{Q}}\big{[}\displaystyle\frac{1}{(s\bar{\tau})^{6}}\big{]}\cdot\bar{I}_{[0,T]}(\psi_{1})\leq \tilde{C}\displaystyle\int_{\omega\times(0,T)}|\psi_{1}|^{2}dxdt.
\end{array}
\end{eqnarray}
Further, by (\ref{z4.36})~and~(\ref{z4.37}), we have
\begin{equation}\label{z4.38}
\displaystyle\sum_{i=1}^{2}\int_{Q}(|\gamma_{1}^{i}|^{2}+|\gamma_{2}^{i}|^{2})dxdt\leq\tilde{ C} \displaystyle\int_{\omega\times(0,T)}|\psi_{1}|^{2}dxdt.
\end{equation}
Combining (\ref{z4.35})~with~(\ref{z4.38}), the desired estimate 
(\ref{z4.5}) follows.
\endpf

\end{document}